\documentclass{article}
 \usepackage{amsmath}

\newtheorem{thm}{Theorem}[section]

\newtheorem{prop}[thm]{Proposition}

\newtheorem{lem}[thm]{Lemma}

\newtheorem{rem}[thm]{Remark}

\newtheorem{ex}[thm]{Example}

\newcommand{\be}{\begin{equation}}
\newcommand{\ee}{\end{equation}}
\newcommand{\ben}{\begin{enumerate}}
\newcommand{\een}{\end{enumerate}}
\newcommand{\beq}{\begin{eqnarray}}
\newcommand{\eeq}{\end{eqnarray}}
\newcommand{\beqn}{\begin{eqnarray*}}
\newcommand{\eeqn}{\end{eqnarray*}}

\newcommand{\pa}{\partial}

\newcommand{\qed}{\hspace*{\fill}Q.E.D.}  %Use at end of proof

\begin{document}
\title{ On Concircular Transformations in Finsler Geometry}
\author{Zhongmin Shen and Guojun Yang\footnote{Corresponding author, supported by the
National Natural Science Foundation of China (11471226) }}
\date{}
\maketitle
\begin{abstract}
 A geodesic circle in Finsler geometry is a
 natural extension of that in a Euclidean
 space. In this paper, we  study geodesic circles and (infinitesimal) concircular
 transformations
 on a Finsler manifold. We characterize a concircular vector
 field with some PDEs on the tangent bundle, and then we obtain respectively necessary and sufficient conditions
 for a concircular vector
 field to be  conformal and a conformal vector field to be
 concircular. We also show conditions for two conformally related
 Finsler metrics to be concircular, and obtain some
 invariant curvature properties under conformal and concircular transformations.


 {\bf MR(2000) subject classification: }
53B40, 53C60
\end{abstract}

\section{Introduction}

A geodesic circle in a Euclidean space is a straight line or a
circle with finite positive radius, and it can be generalized
naturally to Riemann geometry by using Levi-Civita connection
(\cite{NY}), or more generally generalized to Finsler geometry by
the so called Cartan  $Y$-connection introduced by M. Matsumoto in
\cite{AIM} \cite{Ma} (also see Section \ref{sec2} below). A curve
$\gamma=\gamma(s)$ on a Finsler manifold $(M,F)$ with $s$ being
the arc-length is called a geodesic circle if it satisfies
 \be\label{y01}
 D^*_{\dot{\gamma}}D^*_{\dot{\gamma}}\dot{\gamma}
 +g_{\dot{\gamma}}(D^*_{\dot{\gamma}}\dot{\gamma},D^*_{\dot{\gamma}}\dot{\gamma})\
 \dot{\gamma}=0,
 \ee
where $D^*$ is the Cartan $Y$-connection (induced by
$\dot{\gamma}$) and $g_{\dot{\gamma}}$ is the inner product
induced by $F$.
 For two Finsler manifolds
$(M,F)$ and $(\widetilde{M},\widetilde{F})$, a diffeomorphism
$\varphi$ from $(M,F)$ to $(\widetilde{M},\widetilde{F})$ is said
to be  concircular if $\varphi$ maps geodesic circles to  geodesic
circles. For convenience, we say two Finsler metrics on a same
manifold are concircular if they have the same geodesic circles as
points set. Correspondingly,  a vector field $V$ on a Finsler
manifold $(M,F)$ is said to be concircular if   its flow induces
infinitesimal concircular transformations.

A Finsler metric $F=F(x,y)$ with $x\in M, y\in
  T_xM$ defines its fundamental metric
tensor $g_{ij}$ (while $g^{ij}$ the inverse), Cartan torsion
$C_{jk}^i$ and mean Cartan torsion $I_i$ respectively by
 $$g_{ij}:=\dot{\pa}_i\dot{\pa}_j\big(F^2/2\big), \
 \ 2C_{ijk}:=\dot{\pa}_kg_{ij}, \ \ I_i:=g^{jk}C_{ijk}=C^r_{ir}, \
 \ \ \big(\dot{\pa}_i:=\pa/\pa y^i\big).
 $$
It is well known that a Finsler metric is Riemannian iff. the
Cartan torsion, or the mean Cartan torsion vanishes (\cite{De}).
In a Minkowski Finsler space, the geodesic circle equation
(\ref{y01}) is reduced to a simple equation which is closely
related to the Cartan torsion (Example \ref{ex35} below).
 A vector field $V$ on a manifold $M$ induces a flow $\varphi_t$
acting on $M$, and $\varphi_t$ is naturally lifted to a flow
$\widetilde{\varphi}_t$ on the tangent bundle $TM$, where
$\widetilde{\varphi}_t: TM \mapsto TM$ is defined  by
  $\widetilde{\varphi}_t(x,y):=(\varphi_t(x),\varphi_{t*}(y))$. Taking the derivative of $\widetilde{\varphi}_t$ with
  respect to $t$ at $t=0$, we obtain a vector field $V^c$ on the
  tangent bundle $TM$, which is called the complete lift of $V$.
  A vector field $V$ on a Finsler manifold $(M,F)$ is said to be conformal if $F$ keeps
  conformally related under the flow $\widetilde{\varphi}_t$, that is, it
  holds
  $
 \widetilde{\varphi}_t^*F=e^{\sigma_t}F,
 $
where $\sigma_t$ is a function on $M$ for every $t$, and then by
taking the derivative of $\sigma_t$ at $t=0$ we obtain a scalar
function $\rho$ (on $M$) called  a conformal factor. For some
studies on conformal vector fields, one may refer to \cite{MH}
\cite{Y3} \cite{Y4}, for instance.

In 1940s, Yano introduced concircular transformations of
Riemannian manifolds and developed the theory of concircular
geometry in a series of papers (\cite{Yan}). After that, some
researchers did further jobs on concircular transformations in
Riemann geometry (see for instance \cite{Is} \cite{IT} \cite{Ta}
\cite{Vo}). Vogel shows that a concircular transformation of
Riemannian manifolds is a conformal transformation (\cite{Vo}),
and Ishihara proves that a concircular vector field on a
Riemannian manifold is a conformal vector field (\cite{Is}).

For Finsler manifolds, one is wondering whether a concircular
transformation, or a concircular vector field is still conformal.
Some investigations are made in \cite{BS} \cite{JB}. We find that
the Finslerian case is much more complicated since it is closely
related to the Cartan torsion. In this paper, we will first
characterize a concircular vector field by some PDEs (Theorem
\ref{th61} below), and then using Theorem \ref{th61}, we obtain
the following Theorems \ref{th1} and \ref{th2}.

\begin{thm}\label{th1}
 A concircular vector field $V$ on a Finsler manifold is conformal
 if and only if
 the Lie derivative of the mean Cartan torsion along $V^c$ vanishes.
\end{thm}

\begin{thm}\label{th2}
 On a Finsler manifold, a conformal vector field   with the conformal factor $\rho$ is
concircular if and only if $\rho$ satisfies
   \be\label{y1}
 \rho_{i|j}=\lambda g_{ij},\ \ \ \rho^rC_{ri}^k=0, \ \ \ \ \  (\rho_i:=\rho_{x^i},\
 \rho^i:=g^{ir}\rho_r),
  \ee
 where $\lambda=\lambda(x)$ is a scalar function on $M$ and the symbol $_|$
 means the horizontal covariant derivative of Cartan (or  Chern) connection.
\end{thm}

In (\ref{y1}), the horizontal covariant derivative of Cartan
connection can also be replaced by that of Berwald connection due
to the second equation of (\ref{y1}). Theorems \ref{th1} and
\ref{th2} show that a concircular vector field is closely related
to the Cartan torsion or mean Cartan torsion.  For a Riemann
metric, the Cartan torsion and the mean Cartan torsion both
vanish. Then Theorems \ref{th1} and \ref{th2} show that a vector
field $V$ on a Riemann manifold is concircular iff. $V$ is
conformal with the conformal factor $\rho$ satisfying the first
formula in (\ref{y1}) (see \cite{Is}). In Section \ref{sec7}, we
will see  that on certain Finsler manifolds, there are concircular
(resp. conformal) but not conformal (resp. concircular)  vector
fields.

For concircular transformations between two Finsler metrics we
have the following result.

\begin{thm}\label{th3}
 Let $\widetilde{F}$ and $F$ be two conformally related Finsler metrics on a same
 manifold $M$ with $\widetilde{F}=u^{-1}F$. Then we have
 \ben
 \item[{\rm (i)}]  $\widetilde{F}$ and
 $F$ are concircular if and only if
  \be\label{y3}
  u_{i|j}=\lambda g_{ij},\ \ \ u^rC_{ri}^k=0, \ \ \ \ \  (u_i:=u_{x^i},\
 u^i:=g^{ir}u_r),
  \ee
  where $\lambda=\lambda(x)$ is a scalar function on $M$ and the symbol $_|$
 means the horizontal covariant derivative of Cartan (or  Chern) connection of $F$.

 \item[{\rm (ii)}] If $F$ and $\widetilde{F}$ are concircular,
 then $F$ and $\widetilde{F}$ keep the invariance of their
 features of being of scalar (resp. isotropic) flag
 curvature, or of constant flag curvature (in $dim(M)\ge 3$), or an Einstein metric.
 In this case, we have the following formula
  \be\label{y4}
 \widetilde{K}=Ku^2+2\lambda u-u_mu^m,
  \ee
  where $\lambda$ is given by (\ref{y3}), and $K$ (resp. $\widetilde{K}$) denotes the flag curvature or Ricci scalar
  of $F$ (resp. $\widetilde{F}$).
 \een
\end{thm}

Theorem \ref{th3} (i) is  an analogue of Theorem \ref{th2}. In
Theorem \ref{th3} (i), if $F$ is locally Euclidean, then the local
structure of $\widetilde{F}$ can be determined by solving
(\ref{y3}), and this case can be an example to show that a
geodesic (resp. circle) may be mapped to a circle (resp. geodesic)
(see Remark \ref{rem61} following the proof of Theorem \ref{th3}).
Theorem \ref{th3} (ii) provides a similar result as  a projective
map keeps scalar flag curvature unchanged. We are not sure whether
the converse of Theorem \ref{th3} (ii) is true in dimension $n\ge
3$, which holds however at least in Riemnnian case (\cite{Fe}). In
Theorem \ref{th3} (ii), if $F$ is locally Minkowskian,  then
$\widetilde{F}$ is of isotropic flag curvature.

\

We organize the paper as follows. In Section \ref{sec2}, we
introduce the definition of  Cartan $Y$-connection and its basic
properties. In Section \ref{sec3}, we introduce the definition of
 geodesic circles and show some basic properties of the ODE
related to geodesic circles. In Section \ref{sec4}, we show the
notion of Lie derivative and some useful formulas related to Lie
derivative are given. In Section \ref{sec6}, we give the proofs of
our main results, and therein, we also establish the
characterization theorem (Theorem \ref{th61}) for concircular
vector fields. In Section \ref{sec7}, we give some examples
supplementary to Theorems \ref{th1} and \ref{th2}.

\section{Finsler Connections and Cartan $Y$-connection}\label{sec2}

 A spray ${\bf G}$ is a global vector field defined on the
  tangent bundle $TM$,
 $$
 {\bf G}:=y^i\frac{\pa}{\pa x^i}-2G^i\frac{\pa}{\pa y^i},
 $$
where $G^i$ are called the spray (or geodesic) coefficients. Put
 $$
 \delta_i:=\frac{\pa}{\pa x^i}-G^r_i\frac{\pa}{\pa y^r},\ \ \
 \dot{\pa}_i:=\frac{\pa}{\pa y^i},\ \  \   \pa_i:=\frac{\pa}{\pa x^ i},\ \ \ \  \delta
 y^i:=dy^i+G_r^idx^r, \ \ \ \big(G^r_i:=\dot{\pa}_iG^r\big).
 $$
 Then $\{\delta_i,\dot{\pa}_i\}$ is a local frame on the manifold
 $TM$ and $\{dx^i,\delta y^i\}$ is its dual.
We denote by $\pi:\ TM\mapsto M$  the natural projection. Let
$\mathcal{H}$ and $\mathcal{V}$ be two maps from $\pi^*TM$ to
$TTM$ and they are locally given by
 $$
 \mathcal{H}v=v^i\delta_i,\ \ \
 \mathcal{V}v=v^i\dot{\pa}_i,\ \ \  \  \big(v=v^i(x,y)\pa_i\big).
 $$

Let $D$ be a linear connection defined on the pull-back vector
bundle $\pi^*TM$ with $TM$ as the base manifold. We can put
 $$
 D(\frac{\pa}{\pa x^i})=\omega^r_i\frac{\pa}{\pa x^r}=(\Gamma^r_{ik}dx^k+V^r_{ik}\delta  y^k)\frac{\pa}{\pa
 x^r}.
 $$
For a spray tensor $T_i^j$, as an example, the $h$- and
$v$-covariant derivatives (denoted by $_|$ and $|$ respectively)
are defined respectively by
 $$
 T^j_{i|k}:=\delta_kT^j_i+T^r_i\Gamma^j_{rk}-T^j_r\Gamma^r_{ik},\
 \ \ T^j_i|_k:=\dot{\pa}_kT^j_i+T^r_iV^j_{rk}-T^j_rV^r_{ik}.
 $$
The $hh$-curvature ${\bf R}$ and the $hv$-curvature ${\bf P}$ are
given by
 \beqn
{\bf
R}(X,Y)Z:\hspace{-0.6cm}&&=-D_{\mathcal{H}X}D_{\mathcal{H}Y}Z+D_{\mathcal{H}Y}D_{\mathcal{H}X}Z+D_{[\mathcal{H}X,\mathcal{H}Y]}Z,\\
{\bf
P}(X,Y)Z:\hspace{-0.6cm}&&=-D_{\mathcal{H}X}D_{\mathcal{V}Y}Z+D_{\mathcal{V}Y}D_{\mathcal{H}X}Z+D_{[\mathcal{H}X,\mathcal{V}Y]}Z,
 \eeqn
for $X,Y,Z\in \pi^*TM$. Under the natural local basis $\{\pa_i\}$,
we have
 \beqn
R^{\ r}_{k\
ij}\hspace{-0.6cm}&&=\delta_j\Gamma^r_{ki}+\Gamma^s_{ki}\Gamma^r_{sj}-\delta_i\Gamma^r_{kj}-\Gamma^s_{kj}\Gamma^r_{si}
+(\delta_jG^s_i-\delta_iG^s_j)V^r_{ks},\\
P^{\ r}_{k\
ij}\hspace{-0.6cm}&&=\dot{\pa}_j\Gamma^r_{ki}-V^r_{kj|i}-(\Gamma^s_{ji}-G_{ij}^s)V^r_{ks},\ \ \ \ \  (G_{ij}^s:=\dot{\pa}_jG^s_i), \\
 && \big({\bf R}(\pa_i,\pa_j)\pa_k=R^{\ r}_{k\ ij}\pa_r,\ \ \ \   {\bf P}(\pa_i,\pa_j)\pa_k=P^{\ r}_{k\
 ij}\pa_r\big).
 \eeqn

For a Finsler metric $F$, there are three well-known connections:
Cartan, Berwald and Chern connections, which are defined
respectively by putting
 \beq
 &&\Gamma^k_{ij}:=^*\hspace{-0.1cm}\Gamma^k_{ij},\ \
 V^k_{ij}:=C^k_{ij}; \ \ \ \   \Gamma^k_{ij}:=G^k_{ij},\ \
 V^k_{ij}:=0; \ \ \ \   \Gamma^k_{ij}:=^*\hspace{-0.1cm}\Gamma^k_{ij},
  \ \  V^k_{ij}:=0,\nonumber\\
  && \big(\
  ^*\Gamma^k_{ij}:=\frac{1}{2}g^{kl}(\delta_ig_{jl}+\delta_jg_{il}-\delta_lg_{ij}),\
  \ \ \   G^i:=\frac{1}{4}g^{il}\big
  \{[F^2]_{x^ky^l}y^k-[F^2]_{x^l}\big\}\  \big).\label{y5}
 \eeq

In this paper we use Cartan connection $D$ as a tool and the
symbols $_|$ and $|$ denote its $h$- and $v$-covariant derivatives
respectively. Using Cartan connection, we can define  the so
called Cartan $Y$-connection (see \cite{AIM} \cite{Ma}). Let
$Y=Y^i(x)\partial/\partial x^i$ be a non-zero tangent vector filed
on a domain  of the manifold $M$ and $g^*_{ij}(x):=g_{ij}(x,Y(x))$
be the $Y$-Riemannian metric induced from the vector field $Y$.
The Cartan $Y$-connection (or called Barthel  connection), denoted
by $D^*$, is a linear connection on the tangent bundle $TM$ over
the base manifold $M$, with the connection coefficients given by
 \beq
 \Gamma_{jk}^{*i}(x)\hspace{-0.6cm}&&=^*\hspace{-0.1cm}\Gamma_{jk}^i(x,Y(x))+C^i_{jr}(x,Y(x))Y_k^r(x,Y(x)),\label{y6}\\
 &&\big(Y_j^i(x,y):=Y^i_{|j}(x,y)=(\partial_jY^i)(x)+G_j^i(x,y)\big).\nonumber
 \eeq
 We use the symbol $_/$ to denote the covariant derivative of the Cartan $Y$-connection. For a spray tensor $T_i(x,y)$ (as an example), let
 $T^*_i(x):=T_i(x,Y(x))$. Then $T^*_i$ is considered as a tensor on $M$ and we have
  \beq
&&T^*_{i/j}=\big(T_{i|j}+T_i|_rY^r_j\big)|_{y=Y},\label{y7}\\
 &&T^*_{i|j}=\big\{T_{i|j}+T_{i\cdot r}Y^r_i\}|_{y=Y},\ \ \
 (T_{i\cdot r}:=\dot{\pa}_rT_i),\nonumber\\
 &&T^*_{i|j}+T^*_i|_rY^r_j=\{T_{i|j}+T_i|_rY^r_j\big\}_{y=Y}.\label{y8}
 \eeq
Since the Cartan connection is $F$-metric-compatible
($g_{ij|k}=0$, $g_{ij}|_k=0$), the Cartan $Y$-connection is
$g^*$-metric-compatible ($g^*_{ij/k}=0$) by (\ref{y7}).

For a Finsler manifold $(M,F)$  and a curve $\gamma=\gamma(t)$ on
$M$, we always in this paper let $Y$ be a vector field in the
neighborhood of $\gamma$ which is an extension of
$\dot{\gamma}:=d\gamma/dt$ and let $D^*$ be the Cartan
$Y$-connection related to the vector field $Y$. Let
$\dot{\widetilde{\gamma}}:=(\dot{\gamma},\ddot{\gamma})$ be the
tangent vector of the curve
$\widetilde{\gamma}:=(\gamma,\dot{\gamma})$ on $TM$. Then for a
spray tensor $T=T_idx^i$ we have
 \be\label{y9}
D^*_{\dot{\gamma}}T^*=D_{\dot{\widetilde{\gamma}}}T^*=D_{\dot{\widetilde{\gamma}}}T,
 \ee
which follows from (note that
$\dot{\gamma}^rY^k_r=\ddot{\gamma}^k+2G^k$ and $y$ takes the value
$\dot{\gamma}$)
 \beqn
 (D^*_{\dot{\gamma}}T^*)_i\hspace{-0.6cm}&&=\dot{\gamma}^kT^*_{i/k}\stackrel{(\ref{y7})}{=}
 \dot{\gamma}^k\big(T_{i|k}+T_i|_rY^r_k\big)=\dot{\gamma}^kT_{i|k}+(\ddot{\gamma}^k+2G^k)T_i|_k\\
 &&=\big(D_{\dot{\gamma}^k\delta_k+(\ddot{\gamma}^k+2G^k)\dot{\pa}_k}T\big)_i
 =\big(D_{\dot{\gamma}^k\pa_k+\ddot{\gamma}^k\dot{\pa}_k}T\big)_i=\big(D_{\dot{\widetilde{\gamma}}}T\big)_i,
 \eeqn
and similarly
$(D^*_{\dot{\gamma}}T^*)_i=(D_{\dot{\widetilde{\gamma}}}T^*)_i$
from (\ref{y8}). Let $U$ and $V$ be two vector fields along the
curve $\gamma$, and then we have (since $D^*$ is
$g^*$-metric-compatible)
 \be\label{y10}
 \frac{d}{dt}g_{\dot{\gamma}}(U,V)=D^*_{\dot{\gamma}}g_{\dot{\gamma}}(U,V)=
 g_{\dot{\gamma}}(D^*_{\dot{\gamma}}U,V)+g_{\dot{\gamma}}(U,D^*_{\dot{\gamma}}V).
 \ee

\begin{rem}
 In (\ref{y10}), we can replace $D^*_{\dot{\gamma}}$ by
 $D_{\dot{\widetilde{\gamma}}}$ from (\ref{y9}), but can not by  $D_{\dot{\gamma}}$.
\end{rem}

\section{Geodesic circles}\label{sec3}

By a simple observation, we have the following lemma.

\begin{prop}\label{lem31}
 Let $\gamma=\gamma(s)$ be parameterized by the arc-length $s$
 satisfying the following ODE
  $$
  D^*_{\dot{\gamma}}D^*_{\dot{\gamma}}\dot{\gamma}
 +\tau(s)\dot{\gamma}=0, \ \ \  (\dot{\gamma}:=d\gamma/ds),
 $$
 where $\tau$ is a smooth function along $\gamma$. Then we have
 $\tau=g_{\dot{\gamma}}(D^*_{\dot{\gamma}}\dot{\gamma},D^*_{\dot{\gamma}}\dot{\gamma})$.
\end{prop}

{\it Proof :} By (\ref{y10}), we have
 $$
 \tau=-g_{\dot{\gamma}}(D^*_{\dot{\gamma}}D^*_{\dot{\gamma}}\dot{\gamma},\dot{\gamma})
 =g_{\dot{\gamma}}(D^*_{\dot{\gamma}}\dot{\gamma},D^*_{\dot{\gamma}}\dot{\gamma}),
 $$
where we have used
$g_{\dot{\gamma}}(D^*_{\dot{\gamma}}\dot{\gamma},\dot{\gamma})=0$
following from $g_{\dot{\gamma}}(\dot{\gamma},\dot{\gamma})=1$.
\qed

\

Now consider a curve $\gamma=\gamma(t)$ on a Finsler manifold
$(M,F)$ satisfying the ODE
 \be\label{y11}
 D^*_{\dot{\gamma}}D^*_{\dot{\gamma}}\dot{\gamma}
 +g_{\dot{\gamma}}(D^*_{\dot{\gamma}}\dot{\gamma},D^*_{\dot{\gamma}}\dot{\gamma})\
 \dot{\gamma}=0,\ \ \  (\dot{\gamma}:=d\gamma/dt).
 \ee
For the local expansion of the first term in (\ref{y11}),  we have
 \beq
D^*_{\dot{\gamma}}D^*_{\dot{\gamma}}\dot{\gamma}\hspace{-0.6cm}&&=
D^*_{\dot{\gamma}}\big[(\ddot{\gamma}^i+2G^i)\frac{\pa}{\pa
x^i}\big]\nonumber\\
&&\stackrel{(\ref{y6})}{=}\big[\dddot{\gamma}^i+2(\pa_jG^i)\dot{\gamma}^j+2G^i_j\ddot{\gamma}^j\big]\frac{\pa}{\pa
x^i}+(\ddot{\gamma}^i+2G^i)\dot{\gamma}^j(^*\Gamma^k_{ij}+C^k_{ir}Y^r_j)\frac{\pa}{\pa
x^k}\nonumber\\
&&=\Big\{\dddot{\gamma}^k+2(\pa_jG^k)\dot{\gamma}^j-4G^k_jG^j+(\ddot{\gamma}^i+2G^i)\big[3G^k_i+C^k_{ir}(\ddot{\gamma}^r+2G^r)\big]\Big\}\frac{\pa}{\pa
x^k}.\label{y12}
 \eeq

To prove Theorems \ref{th1}--\ref{th3} and Theorem \ref{th61}, we
need the following Proposition \ref{prop31}.

\begin{prop}\label{prop31}
 Arbitrarily fix two vectors $u,v\in T_xM$ with
 $F(u)=1$ and $g_u(u,v)=0$. There is a unique curve $\gamma=\gamma(t)$
 satisfying the ODE (\ref{y11}) with the initial condition $\gamma(0)=x,\dot{\gamma}(0)=u$
 and $D^*_u\dot{\gamma}=v$. For the unique curve
 $\gamma=\gamma(t)$, we have $F(\dot{\gamma})=1$, that is, $t$ is
 the arc-length parameter.
\end{prop}

{\it Proof :} The ODE (\ref{y11}) is of degree three by
(\ref{y12}), and so the uniqueness is obvious. Put
 $$f(t):=g_{\dot{\gamma}}(\dot{\gamma},\dot{\gamma}).
 $$
Then by (\ref{y11}), it easily follows from (\ref{y10}) that
 \be\label{y13}
 f''(t)=g(t)[1-f(t)], \ \ \  f(0)=1,\ \ f'(0)=0,
 \ee
where
$g(t):=2g_{\dot{\gamma}}(D^*_{\dot{\gamma}}\dot{\gamma},D^*_{\dot{\gamma}}\dot{\gamma})$.
Then by an ODE theory, (\ref{y13}) has a unique solution $f(t)=1$,
which implies that $F(\dot{\gamma})=1$.    \qed

\begin{prop}\label{prop32}
 For the ODE (\ref{y11}), if $D^*_{\dot{\gamma}}\dot{\gamma}=0$,
 then $\gamma$ is a geodesic; if $D^*_{\dot{\gamma}}\dot{\gamma}\ne0$,
 then $F(\dot{\gamma})=1$
 iff. $g_{\dot{\gamma}}(D^*_{\dot{\gamma}}\dot{\gamma},D^*_{\dot{\gamma}}\dot{\gamma})=k^2$
 with $k$ being a positive constant.
\end{prop}

 {\it Proof :}  We only consider the case
 $D^*_{\dot{\gamma}}\dot{\gamma}\ne0$. If $F(\dot{\gamma})=1$,
 then we have
 $g_{\dot{\gamma}}(D^*_{\dot{\gamma}}\dot{\gamma},\dot{\gamma})=0$.
 Then by (\ref{y10}) and (\ref{y11}), we obtain
  $$
D^*_{\dot{\gamma}}g_{\dot{\gamma}}(D^*_{\dot{\gamma}}\dot{\gamma},D^*_{\dot{\gamma}}\dot{\gamma})
=2g_{\dot{\gamma}}(D^*_{\dot{\gamma}}D^*_{\dot{\gamma}}\dot{\gamma},D^*_{\dot{\gamma}}\dot{\gamma})
=-2g_{\dot{\gamma}}(D^*_{\dot{\gamma}}\dot{\gamma},D^*_{\dot{\gamma}}\dot{\gamma})
\cdot
g_{\dot{\gamma}}(\dot{\gamma},D^*_{\dot{\gamma}}\dot{\gamma})=0,
  $$
which implies
$g_{\dot{\gamma}}(D^*_{\dot{\gamma}}\dot{\gamma},D^*_{\dot{\gamma}}\dot{\gamma})=constant$.

Conversely, if
$g_{\dot{\gamma}}(D^*_{\dot{\gamma}}\dot{\gamma},D^*_{\dot{\gamma}}\dot{\gamma})=k^2$
is a non-zero constant, then we have
  $$
0=D^*_{\dot{\gamma}}g_{\dot{\gamma}}(D^*_{\dot{\gamma}}\dot{\gamma},D^*_{\dot{\gamma}}\dot{\gamma})
=2g_{\dot{\gamma}}(D^*_{\dot{\gamma}}D^*_{\dot{\gamma}}\dot{\gamma},D^*_{\dot{\gamma}}\dot{\gamma})
=-2k^2
g_{\dot{\gamma}}(\dot{\gamma},D^*_{\dot{\gamma}}\dot{\gamma}),
  $$
which shows
$g_{\dot{\gamma}}(\dot{\gamma},D^*_{\dot{\gamma}}\dot{\gamma})=0$.
By this fact, further we have
 $$
 k^2=g_{\dot{\gamma}}(D^*_{\dot{\gamma}}\dot{\gamma},D^*_{\dot{\gamma}}\dot{\gamma})
 =D^*_{\dot{\gamma}}g_{\dot{\gamma}}(\dot{\gamma},D^*_{\dot{\gamma}}\dot{\gamma})-
g_{\dot{\gamma}}(\dot{\gamma},D^*_{\dot{\gamma}}D^*_{\dot{\gamma}}\dot{\gamma})
=k^2g_{\dot{\gamma}}(\dot{\gamma},\dot{\gamma}),
 $$
from which we see $g_{\dot{\gamma}}(\dot{\gamma},\dot{\gamma})=1$.
\qed

\begin{rem}
 By Proposition \ref{prop32}, the circles of a Finsler manifold
 are determined by the following ODE with an initial condition:
  \be\label{y14}
 D^*_{\dot{\gamma}}D^*_{\dot{\gamma}}\dot{\gamma}
 +k^2\dot{\gamma}=0,\ \ \big(F(\dot{\gamma})=1,\  {\rm or} \ \
 g_{\dot{\gamma}}(D^*_{\dot{\gamma}}\dot{\gamma},D^*_{\dot{\gamma}}\dot{\gamma})=k^2\big),
  \ee
  where $k>0$ is a constant. The number $1/k$ is called the radius
  of the circle.
\end{rem}

In a Minkowski Finsler space, the circle equation (\ref{y14}) is
reduced to    a relatively simple form which is closely related to
the Cartan torsion.

\begin{ex}\label{ex35}
 Let $(R^n,F)$ be a Minkowski space. Then by  the spray $G^i=0$ and (\ref{y12}),  the
 circle equation (\ref{y14}) becomes
 \be\label{y15}
 \dddot{\gamma}^k+\ddot{\gamma}^i\ddot{\gamma}^rC_{ir}^k(\dot{\gamma})+k^2\dot{\gamma}^k=0,
 \ \ \  F(\dot{\gamma})=1.
 \ee
\end{ex}

\begin{ex}\label{ex36}
 Let $(R^n,F)$ be a Euclidean space with $F=|y|$. Then the Cartan torsion vanishes and the
 circle equation (\ref{y15}) becomes
  $$
  \frac{d^3\gamma^i}{ds^3}+k^2\frac{d\gamma^i}{ds}=0.
  $$
  Solving the above ODE we obtain
   \be\label{y16}
 \gamma^i=a^i\cos ks+b^i\sin ks+c^i,
   \ee
   where $a,b,c$ are constant vectors. It is easy to see that
   $F(\dot{\gamma})=1$
   is equivalent to
    \be\label{y17}
 |a|=|b|=\frac{1}{k},\ \ \ \langle a,b\rangle =0.
    \ee
    So by (\ref{y17}), the curve given by (\ref{y16}) is a
    Euclidean circle in the plane spanned by the vectors $a,b$, with the
    center at the point $c$ and the radius $1/k$.
\end{ex}

Let $\gamma=\gamma(s)$ be a geodesic circle with $s$ being the
arc-length parameter. Then $\gamma$ satisfies the ODE (\ref{y01}).
Now let $\gamma$ be parameterized by a general parameter $t$, and
define $\gamma':=d\gamma/dt$. A simple computations shows
 \beqn
&&\gamma'=F(\gamma')\dot{\gamma},\ \ \ \ \
D^*_{\gamma'}\gamma'=F^2(\gamma')D^*_{\dot{\gamma}}\dot{\gamma}
+\frac{g_{\gamma'}(\gamma',D^*_{\gamma'}\gamma')}{F(\gamma')}\
 \dot{\gamma},\\
 &&D^*_{\gamma'}D^*_{\gamma'}\gamma'=F^3(\gamma')D^*_{\dot{\gamma}}D^*_{\dot{\gamma}}\dot{\gamma}
 +3g_{\gamma'}(\gamma',D^*_{\gamma'}\gamma')D^*_{\dot{\gamma}}\dot{\gamma}
 +\frac{d}{dt}\Big(\frac{g_{\gamma'}(\gamma',D^*_{\gamma'}\gamma')}{F(\gamma')}\Big)\dot{\gamma}.
 \eeqn
Following the above and Proposition \ref{lem31}, we immediately
obtain the following Proposition \ref{prop37}, which  will be used
to prove Theorem \ref{th3}.

\begin{prop}\label{prop37}
 A curve $\gamma=\gamma(t)$ under a general parameter $t$ is a
 geodesic circle iff. the following vector $U=U(t)$ along the curve $\gamma$,
 $$
U:=D^*_{\gamma'}D^*_{\gamma'}\gamma'-3\frac{g_{\gamma'}(\gamma',D^*_{\gamma'}\gamma')}{F^2(\gamma')}\
D^*_{\gamma'}\gamma'
 $$
 is  tangent to the curve  $\gamma$.
\end{prop}

\section{Lie derivatives}\label{sec4}

Consider a geometric object $T$ on $M$ ($T$ is not necessarily a
tensor), which is defined along curves on $M$ with the following
form
 \be\label{y18}
 T=T(c)=\big(T^{i_1\cdots}_{j_1\cdots}(c,\dot{c},\ddot{c},\cdots,
 c^{(m)})\big),\ \ \  (c^{(k)}:=d^kc/dt^k),
 \ee
where $c=c(t)$ is an arbitrary curve  parameterized by a general
parameter $t$. Obviously, the value of $T$ at a point $x\in M$ is
dependent on the derivatives of some degrees for a curve passing
through $x$. If $m=0$, then $T$ is defined along points of $M$,
which is the case for  a tensor $T$.
 The components
$T^{i_1\cdots}_{j_1\cdots}$ are determined by local coordinates.
For a map $f:M\mapsto M$, denote by $f_{\#} T$ the local
expression of $T$ under the local coordinate $\widetilde{x} \
(=f(x))$ in $\widetilde{U} \ (=f(U))$. If $T$ is a tensor on $M$,
then $f_{\#}$ coincides with the common map induced from  the
tangent map.

For a vector field $V$  on $M$, it induces a flow
 $\varphi_t$ acting on $M$. The Lie derivative of $T$ along $V$ is defined by (cf. \cite{TI} \cite{Ya})
 \be\label{y19}
 \mathcal{L}_V(T(c)):=\frac{d}{d\epsilon}|_{\epsilon=0}\big[T(\varphi_{\epsilon}(c))-\varphi_{\epsilon{\#}}(T(c))\big].
 \ee
The Lie derivative $\mathcal{L}_VT$ measures the change of $T$
along the vector field $V$.

If $m=1$ in (\ref{y18}), then $T$ is actually defined along points
on $TM$ by putting $(c,\dot{c})=(x,y)$ due to the arbitrariness of
the curve $c=c(t)$, and the vector field $V$ on $M$ is lifted to
 the vector field $V^c$ on $TM$, where $V$ and $V^c$ are
locally related  by
 $$
V=V^i\pa_i,\ \ \ \ \  V^c=V^i\pa_i+y^r(\pa_r V^i)\dot{\pa}_i.
 $$
Denote by $\varphi^c_t$ the flow of $V^c$ acting on $TM$. Then we
have a similar definition for $\mathcal{L}_{V^c}T$ as that in
(\ref{y19}). So in this case we identify $\mathcal{L}_VT$ with
 $\mathcal{L}_{V^c}T$. If $T$ is spray tensor on $M$, for example,
 $T=(T^i_j(x,y))$, by the definition (\ref{y19}), we easily obtain
  \be\label{y20}
 \mathcal{L}_{V^c}T^i_j=V^c(T^i_j)-T^r_j(\pa_rV^i)+T^i_r(\pa_jV^r)=V^rT^i_{j|r}+V^r_{\ |0}T^i_{j\cdot  r}-T^r_jV^i_{\ |r}+T^i_rV^r_{\ |j},
  \ee
where we have used the contraction $V^r_{\ |0}:=V^r_{\ |i}y^i$. By
the definition of Lie derivative, the following Lemma \ref{lem41}
is easily proved.
 \begin{lem}\label{lem41}
  For $y^i$, $g_{ij}$ and $G^i$, we have
 \beqn
 &&\mathcal{L}_{V^c}y^i=0, \ \ \ \ \
 \mathcal{L}_{V^c}g_{ij}=V_{i|j}+V_{j|i}+2V^r_{\ |0}C_{rij},\ \ \ \
 \
\mathcal{L}_{V^c}(2G^i)=V^i_{\  |0|0}+V^rR^i_{\ r},\\
 &&\ \  \big( R^i_{\ k}:=2\pa_k G^i-y^j\pa_jG^i_k+2G^jG^i_{jk}- G^i_j G^j_k\ (the \ Riemann \ curvature)\big).
 \eeqn
\end{lem}

\begin{lem}\label{lem42}
 For Cartan (or Chern) and Berwald connections, we have
  \beqn
 A^i_{jk}:\hspace{-0.5cm}&&=\mathcal{L}_{V^c}(^*\Gamma^i_{jk})=V^i_{\ |j|k}+V^r_{\ |0}F^{\ i}_{j\ kr}+V^rK^{\ i}_{j\
 kr},\\
B^i_{jk}:\hspace{-0.5cm}&&=\mathcal{L}_{V^c}(G^i_{jk})=V^i_{\
;j;k}+V^r_{\ ;0}G^{\ i}_{j\ kr}+V^rH^{\ i}_{j\
 kr},
  \eeqn
  where $K^{\ i}_{j\ kr}$ and $F^{\ i}_{j\ kr}$ are the $hh$- and $hv$-curvatures
  of Chern connection, $H^{\ i}_{j\ kr}$ and $G^{\ i}_{j\ kr}$ are the $hh$- and $hv$-curvatures
  of Berwald connection (see Section \ref{sec2}), and the symbol $_;$ is the $h$-covariant
  derivative of Berwald connection.
\end{lem}

\begin{lem}\label{lem43}
 Related to $A^i_{jk}$ and $B^i_{jk}$ in Lemma \ref{lem42}, we
 have
 \beqn
 &&(\mathcal{L}_{V^c}T_i)_{|j}-\mathcal{L}_{V^c}(T_{i|j})=T_rA^r_{ij}+T_{i\cdot
 r}A^r_{0j},\ \ \ \ \   \mathcal{L}_{V^c}K^{\ m}_{i\ jk}=A^m_{ij|k}+A^r_{0k}F^{\
m}_{i\ jr}-(j/k),\\
 &&(\mathcal{L}_{V^c}T_i)_{;j}-\mathcal{L}_{V^c}(T_{i;j})=T_rB^r_{ij}+T_{i\cdot
 r}B^r_{0j}, \ \ \ \ \   \mathcal{L}_{V^c}H^{\ m}_{i\ jk}=B^m_{ij;k}+B^r_{0k}G^{\
m}_{i\ jr}-(j/k),
 \eeqn
 where $T=(T_i)$ is a spray tensor (as an example), and
 $T_{ij}-(i/j)$ means $T_{ij}-T_{ji}$.
\end{lem}

It is a little lengthy to prove Lemmas \ref{lem42} and \ref{lem43}
(cf. \cite{JB}). We omit the details here.

 \begin{rem}
Acting on a general geometric object $T=(T_i(x,y))$ (as an
example), $\mathcal{L}_{V^c}\pa_j=\pa_j\mathcal{L}_{V^c}$ or
$\mathcal{L}_{V^c}\dot{\pa}_j=\dot{\pa}_j\mathcal{L}_{V^c}$ iff.
it holds respectively ($\widetilde{x}=\varphi_{\epsilon}(x)$)
 $$
 \pa_m\Big\{\big[\varphi_{\epsilon\#}(T_i)-T_r\frac{\pa x^r}{\pa
 \widetilde{x}^i}\big]|_{\epsilon=0}\Big\}\cdot\frac{\pa V^m}{\pa
 x^j}=0,\ \ \
\dot{\pa}_m\Big\{\big[\varphi_{\epsilon\#}(T_i)-T_r\frac{\pa
x^r}{\pa
 \widetilde{x}^i}\big]|_{\epsilon=0}\Big\}\cdot\frac{\pa V^m}{\pa
 x^j}=0.
 $$
If $T$ is a spray tensor or $T$ is the spray $G^i$, the above
conditions are satisfied, because for a spray tensor $T_i$ and the
spray $G^i$ we respectively have
 $$
\varphi_{\epsilon\#}(T_i)-T_r\frac{\pa x^r}{\pa
\widetilde{x}^i}=0,\ \ \ \
\big[\varphi_{\epsilon\#}(G^i)-G^r\frac{\pa \widetilde{x}^i}{\pa
x^r}\big]_{\epsilon=0}=\big[-\frac{1}{2}\frac{\pa^2\widetilde{x}^i}{\pa
x^r\pa x^m}y^ry^m\big]_{\epsilon=0}=0.
 $$
 \end{rem}

\

For a curve $c=c(t)$ with a general parameter $t$, by the
definition (\ref{y19}), we have
 \be\label{y021}
 \mathcal{L}_Vc'=\frac{d}{d\epsilon}(\widetilde{c}'-\varphi_{\epsilon\#}c')
 =\frac{d}{d\epsilon}(\widetilde{c}'-\widetilde{c}')=0,\ \
 \mathcal{L}_Vc''=0, \  \cdots
 \ee
Now consider a geometric object $T$ on $M$ defined along curves,
in the following form
 \be\label{y21}
 T=T(c)=\big(T^{i_1\cdots}_{j_1\cdots}(c,\dot{c},\ddot{c},\cdots,
 c^{(m)})\big),\ \ \  (c^{(k)}:=d^kc/ds^k),
 \ee
where $s$ is the arc-length parameter. When we consider the
covariant derivative of $T$ or $\mathcal{L}_VT$ defined along
curves, the understanding is to regard $T$ or $\mathcal{L}_VT$  as
a new object $T^*=T^*(x)$ defined along points in a neighborhood
of the curve $c$. That is, let $Y$ be a vector field in the
neighborhood of $c$ which is an extension of $\dot{c}:=dc/ds$, and
then taking
$c(s)=x,\dot{c}^i(s)=Y^i(x),\ddot{c}^i(s)=Y^r\pa_rY^i,\cdots$, we
obtain a new geometric object $T^*=T^*(x)$ from $T$ in a
neighborhood of the curve $c$. But we should keep in mind that the
Lie derivative  always acts on an object defined along curves. We
will use the following Proposition \ref{prop45} to prove Theorems
\ref{th61} below.

\begin{prop}\label{prop45}
Let $V$ be a vector field on $M$. For a geometric object $T=(T_i)$
on $M$ defined by (\ref{y21}), we have the following exchanging
formulas,
  \beq
 \dot{c}^k\mathcal{L}_V(T_{j/k})\hspace{-0.6cm}&&=\dot{c}^k(\mathcal{L}_VT_{j})_{/k}-\dot{c}^kT_rA^{*r}_{jk},\label{y22}\
 \ \ \   \big(A^{*r}_{jk}:=\mathcal{L}_V(\Gamma^{*r}_{jk})\big),\\
 \dot{c}^k\mathcal{L}_V(T_{j|k})\hspace{-0.6cm}&&=\dot{c}^k(\mathcal{L}_VT_{j})_{|k}-\dot{c}^kT_rA^{r}_{jk}.\label{y23}
  \eeq
\end{prop}

{\it Proof :}  Note that Lemma \ref{lem43} is of no help in this
proof, and the Lie derivative and all values are taken along the
curve $c$. We only prove (\ref{y22}) for the Cartan $Y$-connection
(it is similar for (\ref{y23})). Let $t$ be a general parameter
  of $c$ with $c':=dc/dt$, and we have
 \be\label{y025}
 \mathcal{L}_V\dot{c}^k=\mathcal{L}_V(F^{-1}(c'^k)\
 c'^k)=(\mathcal{L}_{V^c}F^{-1})c'^k=-(\mathcal{L}_{V^c}\ln F)\dot{c}^k,
 \ee
where we have used $\mathcal{L}_Vc'^k=0$  by (\ref{y021}); or in
another way we have
 $$
 \mathcal{L}_V\dot{c}^k=\mathcal{L}_{V^c}(F^{-1}(y)y^k)
 =(\mathcal{L}_{V^c}F^{-1})y^k=-(\mathcal{L}_{V^c}\ln F)\dot{c}^k.
 $$
Now by $T_{j/k}=\pa_kT_j-T_r\Gamma^{*r}_{jk}$ we have
 \be\label{y24}
\dot{c}^k\mathcal{L}_V(T_{j/k})=\dot{c}^k\mathcal{L}_V(\pa_kT_j)-\dot{c}^k\mathcal{L}_V(T_r\Gamma^{*r}_{jk}).
 \ee
In the right hand side of (\ref{y24}),  the first term is written
as
 \beqn
\dot{c}^k\mathcal{L}_V(\pa_kT_j)\hspace{-0.6cm}&&=\mathcal{L}_V(\dot{c}^k\pa_kT_j)-(\pa_kT_j)(\mathcal{L}_V\dot{c}^k)
\stackrel{(\ref{y025})}{=}\mathcal{L}_V(\frac{d}{ds}T_j)+(\mathcal{L}_{V^c}\ln F)\frac{d}{ds}T_j\\
&&=\frac{d}{ds}(\mathcal{L}_VT_j),
 \eeqn
the last equality of which follows from
 \beqn
\mathcal{L}_V(\frac{d}{ds}T_j)\hspace{-0.6cm}&&=\mathcal{L}_V(F^{-1}\frac{d}{dt}T_j)
=(\mathcal{L}_{V^c}F^{-1})\frac{d}{dt}T_j+F^{-1}(\mathcal{L}_V\frac{d}{dt}T_j)\\
&&=-(\mathcal{L}_{V^c}\ln F)\frac{d}{ds}T_j+F^{-1}\frac{d}{dt}\mathcal{L}_VT_j\\
&&=-(\mathcal{L}_{V^c}\ln
F)\frac{d}{ds}T_j+\frac{d}{ds}\mathcal{L}_VT_j.
 \eeqn
Thus (\ref{y24}) gives
 \be\label{y25}
\dot{c}^k\mathcal{L}_V(T_{j/k})=\frac{d}{ds}(\mathcal{L}_VT_j)-\dot{c}^k\mathcal{L}_V(T_r\Gamma^{*r}_{jk}).
 \ee
On the other hand we have
 \be\label{y26}
\dot{c}^k(\mathcal{L}_VT_j)_{/k}=\dot{c}^k\pa_k(\mathcal{L}_VT_j)-(\mathcal{L}_VT_r)\dot{c}^k\Gamma^{*r}_{jk}
=\frac{d}{ds}(\mathcal{L}_VT_j)-(\mathcal{L}_VT_r)\dot{c}^k\Gamma^{*r}_{jk}
 \ee
Then $(\ref{y25})-(\ref{y26})$ gives
 \beqn
-\dot{c}^k\mathcal{L}_V(T_r\Gamma^{*r}_{jk})+(\mathcal{L}_VT_r)\dot{c}^k\Gamma^{*r}_{jk}
\hspace{-0.6cm}&&=-\dot{c}^k\big[(\mathcal{L}_VT_r)\Gamma^{*r}_{jk}+T_r(\mathcal{L}_V\Gamma^{*r}_{jk})\big]+(\mathcal{L}_VT_r)\dot{c}^k\Gamma^{*r}_{jk}\\
&&=-\dot{c}^kT_r(\mathcal{L}_V\Gamma^{*r}_{jk}).
 \eeqn
This gives the proof of (\ref{y22}).     \qed

\

Using Lie derivative, we can characterize conformal vector fields
and concircular vector fields as follows.

  A vector field $V$ is conformal iff.
 it satisfies
  \be\label{y27}
 \mathcal{L}_{V^c}g_{ij}=2\rho g_{ij},\ \ \  (\rho=\rho(x)\ on \
 M).
  \ee
The scalar function $\rho$ is just the conformal factor. $V$ is
homothetic iff. $\rho$ is a constant. $V$ is Killing iff.
$\rho=0$.

By the meaning of Lie derivative and the definition of concircular
vector fields, we see that  a concircular vector field $V$ is
characterized by the following equation
 \be\label{y29}
   \mathcal{L}_VU=0,
 \ \ \  (U:=D^*_{\dot{\gamma}}D^*_{\dot{\gamma}}\dot{\gamma}
 +g_{\dot{\gamma}}(D^*_{\dot{\gamma}}\dot{\gamma},D^*_{\dot{\gamma}}\dot{\gamma})\
 \dot{\gamma}=0,\ \ \  \dot{\gamma}=d\gamma/ds),
 \ee
where $\gamma=\gamma(s)$ is an arbitrary curve  with $s$ being the
arc-length parameter. We will compute $\mathcal{L}_VU$ in the next
section in the proof of
 Theorem \ref{th61} below.

\section{Proofs of main results}\label{sec6}

\subsection{Characterization of concircular vector fields}

Before we prove Theorems \ref{th1} and \ref{th2}, we first give a
characterization for concircular vector fields by some PDEs.
Define some spray tensors on $TM$ as follows:
 \beq
 T^k_{ij}:\hspace{-0.7cm}&&=3\big[(\mathcal{L}_{V^c}y_i)\delta^k_j+(i/j)\big]
 -2F^2\mathcal{L}_{V^c}C^k_{ij}-2(\mathcal{L}_{V^c}y_r)C^r_{ij}y^k,\label{y30}\\
 \theta^k_i:\hspace{-0.7cm}&&=3\big[F^2(\mathcal{L}_{V^c}y_i)-V^c(F^2)y_i\big]y^k
 +3F^2V^c(F^2)\delta^k_i,\label{y31}\\
 S^k:\hspace{-0.7cm}&&=-[V^c(F^2)]_{|0|0}y^k+2F^2A^k_{00|0},\label{y32}\\
Z^k_i:\hspace{-0.7cm}&&=\big\{4A^r_{00}g_{ir}-[V^c(F^2)]_{|i}-4(\mathcal{L}_{V^c}y_i)_{|0}\big\}y^k
-3[V^c(F^2)]_{|0}\delta^k_i\nonumber\\
&&+2F^2(3A^k_{i0}+2A^r_{00}C^k_{ir}),\label{y33}\\
\lambda^k:\hspace{-0.7cm}&&=4\big[A^r_{00}y_r-2[V^c(F^2)]_{|0}\big]y^k+6F^2A^k_{00}.\label{y34}
 \eeq

\begin{thm}\label{th61}
 Let $V$ be a vector field on a Finsler manifold $(M,F)$. Then $V$
 is concircular iff. $V$ satisfies the following PDEs on the
 tangent bundle $TM$:
  \beq\label{y35}
F^4T^k_{ij}=y_i\theta^k_j+y_j\theta^k_i,\ \ \  \ \ \   S^k=0,\ \ \
\ \ \   F^2Z^k_i=\lambda^ky_i,
  \eeq
  where $T^k_{ij},\theta^k_i,S^k,Z^k_i$ and $\lambda^k$ are given
  by (\ref{y30})--(\ref{y34}).
\end{thm}

Let $\gamma=\gamma(s)$ be an arbitrary curve with $s$ being the
arc-length parameter. We will use the Cartan $Y$-connection as a
tool, where $Y$ is a vector field as an extension of
$\dot{\gamma}$ in a neighborhood of $\gamma$. Since a concircular
vector field $V$ is characterized by (\ref{y29}), to prove Theorem
\ref{th61},  we need to first compute $\mathcal{L}_VU$. Note that
in the following Lemmas \ref{lem62}--\ref{lem64}, all quantities
take values along the curve $\gamma$. For example, we have
$\mathcal{L}_VC^k_{ij}=\mathcal{L}_V\big[C^k_{ij}(\gamma(s),\dot{\gamma}(s))\big]\ne\mathcal{L}_{V^c}C^k_{ij}$,
but  we have $\mathcal{L}_Vg_{ij}=\mathcal{L}_{V^c}g_{ij}$ due to
the zero-homogeneity of $g_{ij}$.

\begin{lem}\label{lem62}
 For  two terms in $\mathcal{L}_VU$ we have
  \beq
 \mathcal{L}_V(D^*_{\dot{\gamma}}\dot{\gamma})^k\hspace{-0.6cm}&&=
 -2(\mathcal{L}_{V^c}\ln F)(D^*_{\dot{\gamma}}\dot{\gamma})^k-\frac{d}{ds}(\mathcal{L}_{V^c}\ln F)\dot{\gamma}^k
 +A^k_{rm}\dot{\gamma}^r\dot{\gamma}^m,\label{y36}\\
\mathcal{L}_V(D^*_{\dot{\gamma}}D^*_{\dot{\gamma}}\dot{\gamma})^k\hspace{-0.6cm}&&=
 -3(\mathcal{L}_{V^c}\ln F)(D^*_{\dot{\gamma}}D^*_{\dot{\gamma}}\dot{\gamma})^k-3\frac{d}{ds}(\mathcal{L}_{V^c}\ln  F)(D^*_{\dot{\gamma}}\dot{\gamma})^k
 -\frac{d^2}{ds^2}(\mathcal{L}_{V^c}\ln F)\dot{\gamma}^k\nonumber\\
 &&+3A^k_{rm}(D^*_{\dot{\gamma}}\dot{\gamma})^r\dot{\gamma}^m+A^k_{rm|i}\dot{\gamma}^r\dot{\gamma}^m\dot{\gamma}^i
 +2A^p_{rm}C^k_{pi}(D^*_{\dot{\gamma}}\dot{\gamma})^i\dot{\gamma}^r\dot{\gamma}^m\nonumber\\
 &&-(\mathcal{L}_{V^c}\ln  F)C^k_{rm}(D^*_{\dot{\gamma}}\dot{\gamma})^r(D^*_{\dot{\gamma}}\dot{\gamma})^m
 +(\mathcal{L}_VC^k_{rm})(D^*_{\dot{\gamma}}\dot{\gamma})^r(D^*_{\dot{\gamma}}\dot{\gamma})^m.\label{y37}
  \eeq
\end{lem}

{\it Proof :} First note that
 $$
\dot{\gamma}^rY^k_r=(D^*_{\dot{\gamma}}\dot{\gamma})^k,\ \ \
 \mathcal{L}_V\dot{\gamma}^k=-(\mathcal{L}_{V^c}\ln F)\dot{\gamma}^k,\ \ \
\dot{\gamma}^m\mathcal{L}_V(C^k_{mr}Y^r_i)=0\ (by\
\dot{\gamma}^mC^k_{mr}=0).
 $$
 By  (\ref{y22}) in Proposition \ref{prop45}, we have
\beqn
\mathcal{L}_V(D^*_{\dot{\gamma}}\dot{\gamma})^k\hspace{-0.6cm}&&=\mathcal{L}_V(\dot{\gamma}^i\dot{\gamma}^k_{/i})
=(\mathcal{L}_V\dot{\gamma}^i)\dot{\gamma}^k_{/i}+\dot{\gamma}^i\mathcal{L}_V(\dot{\gamma}^k_{/i})\\
&&=-(\mathcal{L}_{V^c}\ln
F)\dot{\gamma}^i\dot{\gamma}^k_{/i}+\dot{\gamma}^i(\mathcal{L}_V\dot{\gamma}^k)_{/i}
+\dot{\gamma}^i\dot{\gamma}^m\big[A^k_{mi}+\mathcal{L}_V(C^k_{mr}Y^r_i)\big],
 \eeqn
which immediately gives (\ref{y36}).  To prove (\ref{y37}), first
we have
 \beq
\mathcal{L}_V(D^*_{\dot{\gamma}}D^*_{\dot{\gamma}}\dot{\gamma})^k
\hspace{-0.6cm}&&=\mathcal{L}_V\big[\dot{\gamma}^i(D^*_{\dot{\gamma}}\dot{\gamma})^k_{/i}\big]
=(\mathcal{L}_V\dot{\gamma}^i)(D^*_{\dot{\gamma}}\dot{\gamma})^k_{/i}+
\dot{\gamma}^i\mathcal{L}_V\big[(D^*_{\dot{\gamma}}\dot{\gamma})^k_{/i}\big]\nonumber\\
&&=-(\mathcal{L}_{V^c}\ln
F)(D^*_{\dot{\gamma}}D^*_{\dot{\gamma}}\dot{\gamma})^k+
\dot{\gamma}^i\mathcal{L}_V\big[(D^*_{\dot{\gamma}}\dot{\gamma})^k_{/i}\big].\label{y38}
 \eeq
By (\ref{y22}) and then by (\ref{y36}), the second term in the
right hand side of (\ref{y38}) is given by
 \beq
\dot{\gamma}^i\mathcal{L}_V\big[(D^*_{\dot{\gamma}}\dot{\gamma})^k_{/i}\big]\hspace{-0.6cm}&&=
\dot{\gamma}^i\big[\mathcal{L}_V(D^*_{\dot{\gamma}}\dot{\gamma})^k\big]_{/i}+
\dot{\gamma}^i(D^*_{\dot{\gamma}}\dot{\gamma})^m\big[A^k_{mi}+\mathcal{L}_V(C^k_{mr}Y^r_i)\big]\nonumber\\
&&=-2(\mathcal{L}_{V^c}\ln
F)(D^*_{\dot{\gamma}}D^*_{\dot{\gamma}}\dot{\gamma})^k-3\frac{d}{ds}(\mathcal{L}_{V^c}\ln
F)(D^*_{\dot{\gamma}}\dot{\gamma})^k-\frac{d^2}{ds^2}(\mathcal{L}_{V^c}\ln F)\dot{\gamma}^k\nonumber\\
&&+3A^k_{rm}(D^*_{\dot{\gamma}}\dot{\gamma})^r\dot{\gamma}^m+A^k_{rm/i}\dot{\gamma}^r\dot{\gamma}^m\dot{\gamma}^i
+(D^*_{\dot{\gamma}}\dot{\gamma})^m\dot{\gamma}^i\mathcal{L}_V(C^k_{mr}Y^r_i).\label{y39}
 \eeq
For the last two terms in (\ref{y39}) we have
 \beq
A^k_{rm/i}\dot{\gamma}^r\dot{\gamma}^m\dot{\gamma}^i
\hspace{-0.6cm}&&=A^k_{rm|i}\dot{\gamma}^r\dot{\gamma}^m\dot{\gamma}^i+
A^k_{rm}|_pY^p_i\dot{\gamma}^r\dot{\gamma}^m\dot{\gamma}^i\nonumber\\
&&
=A^k_{rm|i}\dot{\gamma}^r\dot{\gamma}^m\dot{\gamma}^i+A^p_{rm}C^k_{pi}\dot{\gamma}^r\dot{\gamma}^m(D^*_{\dot{\gamma}}\dot{\gamma})^i,\label{y40}\\
\dot{\gamma}^i\mathcal{L}_V(C^k_{mr}Y^r_i)
 \hspace{-0.6cm}&&=\mathcal{L}_V\big((D^*_{\dot{\gamma}}\dot{\gamma})^rC^k_{mr}\big)+(\mathcal{L}_{V^c}\ln  F)C^k_{mr}(D^*_{\dot{\gamma}}\dot{\gamma})^r.\label{y41}
 \eeq
By (\ref{y36}) we have
 \beq
\mathcal{L}_V\big((D^*_{\dot{\gamma}}\dot{\gamma})^rC^k_{mr}\big)
\hspace{-0.6cm}&&=(D^*_{\dot{\gamma}}\dot{\gamma})^r\mathcal{L}_VC^k_{mr}
+C^k_{mr}\mathcal{L}_V(D^*_{\dot{\gamma}}\dot{\gamma})^r\nonumber\\
&&=(D^*_{\dot{\gamma}}\dot{\gamma})^r\mathcal{L}_VC^k_{mr}
+C^k_{mr}\big[-2(\mathcal{L}_{V^c}\ln
F)(D^*_{\dot{\gamma}}\dot{\gamma})^r
 +A^r_{ip}\dot{\gamma}^i\dot{\gamma}^p\big].\label{y42}
 \eeq
Now plugging (\ref{y42}) into (\ref{y41}), then (\ref{y40}) and
(\ref{y41}) into (\ref{y39}), and then (\ref{y39}) into
(\ref{y38}), we finally obtain (\ref{y37}).   \qed

\begin{lem}
 The equation $\mathcal{L}_VU=0$ in (\ref{y29}) is equivalent to
 (under the condition $U=0$)
  \beq
 0\hspace{-0.6cm}&&=\Big\{\big[\mathcal{L}_Vg_{ij}-2(\mathcal{L}_{V^c}\ln  F)g_{ij}\big](D^*_{\dot{\gamma}}\dot{\gamma})^i(D^*_{\dot{\gamma}}\dot{\gamma})^j-\frac{d^2}{ds^2}(\mathcal{L}_{V^c}\ln  F)+2g_{ij}A^i_{rm}\dot{\gamma}^r\dot{\gamma}^m(D^*_{\dot{\gamma}}\dot{\gamma})^j\Big\}\dot{\gamma}^k\nonumber\\
 &&-3\frac{d}{ds}(\mathcal{L}_{V^c}\ln  F)(D^*_{\dot{\gamma}}\dot{\gamma})^k
 +3A^k_{rm}(D^*_{\dot{\gamma}}\dot{\gamma})^r\dot{\gamma}^m+A^k_{rm|i}\dot{\gamma}^r\dot{\gamma}^m\dot{\gamma}^i
 +2A^p_{rm}C^k_{pi}(D^*_{\dot{\gamma}}\dot{\gamma})^i\dot{\gamma}^r\dot{\gamma}^m\nonumber\\
 &&-(\mathcal{L}_{V^c}\ln  F)C^k_{rm}(D^*_{\dot{\gamma}}\dot{\gamma})^r(D^*_{\dot{\gamma}}\dot{\gamma})^m
 +(\mathcal{L}_VC^k_{rm})(D^*_{\dot{\gamma}}\dot{\gamma})^r(D^*_{\dot{\gamma}}\dot{\gamma})^m.\label{y43}
  \eeq
\end{lem}

{\it Proof :} First by the definition of $U$ we have
 \beq
\mathcal{L}_VU^k\hspace{-0.6cm}&&=\mathcal{L}_V(D^*_{\dot{\gamma}}D^*_{\dot{\gamma}}\dot{\gamma})^k+
\big\{(\mathcal{L}_Vg_{ij})(D^*_{\dot{\gamma}}\dot{\gamma})^i(D^*_{\dot{\gamma}}\dot{\gamma})^j
+2g_{ij}(D^*_{\dot{\gamma}}\dot{\gamma})^j\mathcal{L}_V(D^*_{\dot{\gamma}}\dot{\gamma})^i\nonumber\\
&&-(\mathcal{L}_{V^c}\ln
F)g_{\dot{\gamma}}(D^*_{\dot{\gamma}}\dot{\gamma},D^*_{\dot{\gamma}}\dot{\gamma})\big\}
 \dot{\gamma}^k.\label{y44}
 \eeq
Plugging (\ref{y36}), (\ref{y37}) and
 $$
(D^*_{\dot{\gamma}}D^*_{\dot{\gamma}}\dot{\gamma})^k=-g_{\dot{\gamma}}(D^*_{\dot{\gamma}}\dot{\gamma},D^*_{\dot{\gamma}}\dot{\gamma})
 \dot{\gamma}^k
 $$
into (\ref{y44}), we immediately obtain (\ref{y43}) from
$\mathcal{L}_VU=0$.   \qed

\

To simplify (\ref{y43}),  we rewrite (\ref{y43}) in a different
form in the following lemma.

\begin{lem}\label{lem64}
Put $X:=D^*_{\dot{\gamma}}\dot{\gamma}$. Then (\ref{y43}) is
equivalent to
 \be\label{y45}
 \widetilde{T}^k_{ij}X^iX^j+\widetilde{Z}^k_iX^i+\widetilde{S}^k=0,
 \ee
 where $\widetilde{T}^k_{ij},\widetilde{Z}^k_i$ and
 $\widetilde{S}^k$ are defined by
 \beqn
\widetilde{T}^k_{ij}:\hspace{-0.6cm}&&=-\frac{3}{2}\big[(\mathcal{L}_Vg_{ir})\dot{\gamma}^r\delta^k_j+(i/j)\big]
+\mathcal{L}_VC^k_{ij}-(\mathcal{L}_{V^c}\ln F)C^k_{ij}+(\mathcal{L}_Vg_{rm})C^r_{ij}\dot{\gamma}^m\dot{\gamma}^k,\\
\widetilde{Z}^k_i:\hspace{-0.6cm}&&=\big[2A^p_{rm}g_{ip}-\frac{1}{2}(\mathcal{L}_Vg_{rm})_{|i}
-2(\mathcal{L}_Vg_{ir})_{|m}\big]\dot{\gamma}^m\dot{\gamma}^r\dot{\gamma}^k
-\frac{3}{2}(\mathcal{L}_Vg_{rm})_{|p}\dot{\gamma}^m\dot{\gamma}^r\dot{\gamma}^p\delta^k_i\\
&&3A^k_{ir}\dot{\gamma}^r+2A^p_{rm}C^k_{ip}\dot{\gamma}^m\dot{\gamma}^r,\\
\widetilde{S}^k:\hspace{-0.6cm}&&=\big[-\frac{1}{2}(\mathcal{L}_Vg_{ij})_{|r|m}\dot{\gamma}^j\dot{\gamma}^k
+A^k_{rm|i}\big]\dot{\gamma}^m\dot{\gamma}^r\dot{\gamma}^i.
 \eeqn
\end{lem}

{\it Proof :} It needs to expand the derivatives of
$\mathcal{L}_{V^c}\ln
F=(\mathcal{L}_Vg_{ij})\dot{\gamma}^i\dot{\gamma}^j/2$ with
respect to $s$  in (\ref{y43}). We have the following direct
results:
 \beq
 \frac{d}{ds}(\mathcal{L}_{V^c}\ln  F)\hspace{-0.6cm}&&=\frac{1}{2}(\mathcal{L}_Vg_{ij})_{|r}\dot{\gamma}^i\dot{\gamma}^j\dot{\gamma}^r
 +(\mathcal{L}_Vg_{ij})\dot{\gamma}^i(D^*_{\dot{\gamma}}{\dot{\gamma}})^j,\label{y46}\\
 \frac{d^2}{ds^2}(\mathcal{L}_{V^c}\ln F)\hspace{-0.6cm}&&=[\mathcal{L}_Vg_{ij}-2(\mathcal{L}_{V^c}\ln  F)g_{ij}](D^*_{\dot{\gamma}}\dot{\gamma})^i(D^*_{\dot{\gamma}}\dot{\gamma})^j
 +\frac{1}{2}(\mathcal{L}_Vg_{ij})_{|r|m}\dot{\gamma}^i\dot{\gamma}^j\dot{\gamma}^r\dot{\gamma}^m\nonumber\\
 &&+\frac{1}{2}(\mathcal{L}_Vg_{ij})_{|r}\dot{\gamma}^i\dot{\gamma}^j(D^*_{\dot{\gamma}}\dot{\gamma})^r
 +2(\mathcal{L}_Vg_{ij})_{|r}\dot{\gamma}^i\dot{\gamma}^r(D^*_{\dot{\gamma}}\dot{\gamma})^j\nonumber\\
 &&-(\mathcal{L}_Vg_{ij})C^i_{rm}\dot{\gamma}^j(D^*_{\dot{\gamma}}\dot{\gamma})^r(D^*_{\dot{\gamma}}\dot{\gamma})^m.\label{y47}
 \eeq
 Plugging (\ref{y46}) and (\ref{y47}) into (\ref{y43}), we immediately
 obtain (\ref{y45}).
 We can conclude (\ref{y46}) and (\ref{y47}) in the following way. First it is easy to see that
 \be\label{y48}
 (\mathcal{L}_Vg_{ij})|_my^iy^j=0,\ \ \ \ \  (\mathcal{L}_Vg_{ij})_{|r}|_my^iy^jy^r
 =(\mathcal{L}_Vg_{ij})|_{m|r}y^iy^jy^r=0.
 \ee
 Using Cartan $Y$-connection, we see
 $$
 \frac{d}{ds}(\mathcal{L}_{V^c}\ln  F)=\frac{1}{2}(\mathcal{L}_Vg_{ij})_{/r}\dot{\gamma}^i\dot{\gamma}^j\dot{\gamma}^r
 +(\mathcal{L}_Vg_{ij})\dot{\gamma}^i(D^*_{\dot{\gamma}}{\dot{\gamma}})^j,
 $$
which gives (\ref{y46}) from (\ref{y7}) and  (\ref{y48}). To show
(\ref{y47}), we first have
 \beq
&&\ \ \ \  -\frac{d^2}{ds^2}(\mathcal{L}_{V^c}\ln
F)+[\mathcal{L}_Vg_{ij}-2(\mathcal{L}_{V^c}\ln  F)g_{ij}](D^*_{\dot{\gamma}}\dot{\gamma})^i(D^*_{\dot{\gamma}}\dot{\gamma})^j\label{y49}\\
 &&=-\frac{1}{2}(\mathcal{L}_Vg_{ij})_{/r/m}\dot{\gamma}^i\dot{\gamma}^j\dot{\gamma}^r\dot{\gamma}^m
 -\frac{1}{2}(\mathcal{L}_Vg_{ij})_{/r}\dot{\gamma}^i\dot{\gamma}^j(D^*_{\dot{\gamma}}\dot{\gamma})^r
 -2(\mathcal{L}_Vg_{ij})_{/r}\dot{\gamma}^i\dot{\gamma}^r(D^*_{\dot{\gamma}}\dot{\gamma})^j.\nonumber
 \eeq
For the three terms in the right hand side of (\ref{y49}), we
rewrite them as follows. By (\ref{y7}) and (\ref{y48}), we easily
get
 $$
(\mathcal{L}_Vg_{ij})_{/r}\dot{\gamma}^i\dot{\gamma}^j=(\mathcal{L}_Vg_{ij})_{|r}\dot{\gamma}^i\dot{\gamma}^j.
 $$
By (\ref{y7}) and
 \be\label{y052}
 (\mathcal{L}_Vg_{ij})|_m\dot{\gamma}^i
 =\big[2\mathcal{L}_VC_{ijm}-(\mathcal{L}_Vg_{rj})C^r_{im}-(\mathcal{L}_Vg_{ir})C^r_{jm}\big]\dot{\gamma}^i
 =-(\mathcal{L}_Vg_{ir})C^r_{jm}\dot{\gamma}^i,
 \ee
we have
 \beqn
(\mathcal{L}_Vg_{ij})_{/r}\dot{\gamma}^i\dot{\gamma}^r(D^*_{\dot{\gamma}}\dot{\gamma})^j
\hspace{-0.6cm}&&=(\mathcal{L}_Vg_{ij})_{|r}\dot{\gamma}^i\dot{\gamma}^r(D^*_{\dot{\gamma}}\dot{\gamma})^j
+(\mathcal{L}_Vg_{ij})|_mY^m_r\dot{\gamma}^i\dot{\gamma}^r(D^*_{\dot{\gamma}}\dot{\gamma})^j\\
&&=(\mathcal{L}_Vg_{ij})_{|r}\dot{\gamma}^i\dot{\gamma}^r(D^*_{\dot{\gamma}}\dot{\gamma})^j
+(\mathcal{L}_Vg_{ij})|_m\dot{\gamma}^i(D^*_{\dot{\gamma}}\dot{\gamma})^m(D^*_{\dot{\gamma}}\dot{\gamma})^j\\
&&=(\mathcal{L}_Vg_{ij})_{|r}\dot{\gamma}^i\dot{\gamma}^r(D^*_{\dot{\gamma}}\dot{\gamma})^j
-(\mathcal{L}_Vg_{ir})C^r_{jm}\dot{\gamma}^i(D^*_{\dot{\gamma}}\dot{\gamma})^m(D^*_{\dot{\gamma}}\dot{\gamma})^j.
 \eeqn
Finally for the first term of (\ref{y49}), we have
 \beqn
(\mathcal{L}_Vg_{ij})_{/r/m}\dot{\gamma}^i\dot{\gamma}^j\dot{\gamma}^r\dot{\gamma}^m
\hspace{-0.6cm}&&=\big[(\mathcal{L}_Vg_{ij})_{|r}+(\mathcal{L}_Vg_{ij})|_pY^p_r\big]_{/m}\dot{\gamma}^i\dot{\gamma}^j\dot{\gamma}^r\dot{\gamma}^m\\
&&=\Big\{(\mathcal{L}_Vg_{ij})_{|r|m}+(\mathcal{L}_Vg_{ij})_{|r}|_pY^p_m+(\mathcal{L}_Vg_{ij})|_pY^p_{r/m}
\\
&&\ \ \ \ \  \
+\big[(\mathcal{L}_Vg_{ij})|_{p|m}+(\mathcal{L}_Vg_{ij})|_p|_qY^q_m\big]Y^p_r\Big\}\dot{\gamma}^i\dot{\gamma}^j\dot{\gamma}^r\dot{\gamma}^m\\
&&
\hspace{-0.12cm}\stackrel{(\ref{y48})}{=}\big[(\mathcal{L}_Vg_{ij})_{|r|m}+(\mathcal{L}_Vg_{ij})|_p|_qY^q_mY^p_r\big]
\dot{\gamma}^i\dot{\gamma}^j\dot{\gamma}^r\dot{\gamma}^m\\
&&=(\mathcal{L}_Vg_{ij})_{|r|m}\dot{\gamma}^i\dot{\gamma}^j\dot{\gamma}^r\dot{\gamma}^m
+2(\mathcal{L}_Vg_{ij})C^i_{rm}\dot{\gamma}^j(D^*_{\dot{\gamma}}\dot{\gamma})^r(D^*_{\dot{\gamma}}\dot{\gamma})^m,
 \eeqn
in which, the last equality follows from
 \beqn
(\mathcal{L}_Vg_{ij})|_p|_qy^iy^j\hspace{-0.5cm}&&=\big[(\mathcal{L}_Vg_{ij})|_py^iy^j\big]|_q-2(\mathcal{L}_Vg_{qj})|_py^j
\\
&&\hspace{-0.12cm}\stackrel{(\ref{y48})}{=}-2(\mathcal{L}_Vg_{qj})|_py^j\stackrel{(\ref{y052})}{=}2(\mathcal{L}_Vg_{ij})y^jC^i_{pq}.
 \eeqn
Thus we obtain (\ref{y47}) from (\ref{y49}).   \qed

\begin{lem}\label{lem65}
In an $n$-dimensional inner product space with the metric matrix
$(g_{ij})$, let $a_{ij}v^iv^j=0$ be a quadratic-form equation
holding for arbitrary $v\in U^{\bot}$, where $U^{\bot}$ is an
$(n-1)$-dimensional space perpendicular to a unit vector
$u=(u^i)$. Then we have
 \beq\label{y50}
 &&\hspace{2cm}  a_{ij}=\theta_iu_j+\theta_ju_i, \\
 &&\big(\theta_i:=a_{ir}u^r-\frac{1}{2}a_{rm}u^ru^mu_i,\ \ \
 u_i:=g_{ij}u^j\big).\nonumber
 \eeq
\end{lem}

{\it Proof :}  First we have
 $$
 0=a_{ij}(v^i+\bar{v}^i)(v^j+\bar{v}^j)=2a_{ij}v^i\bar{v}^j, \ \ \
 \  (\forall v\in  U^{\bot},\  \bar{v}\in  U^{\bot}),
 $$
which implies $a_{ij}v^j=\lambda u_i$ for some $\lambda
=\lambda(\bar{v})$. Since $u=(u^i)$ is a unit vector, we easily
get $\lambda=a_{ij}u^iv^j$. Thus $a_{ij}v^j=\lambda u_i$ is
written as
 $$
(a_{ij}-a_{jr}u^ru_i)v^j=0,\ \ \ \  (\forall v\in  U^{\bot}),
 $$
which gives
 \be\label{y51}
 a_{ij}-a_{jr}u^ru_i=\tau_i u_j, \ \ \  (for\ some\
 \tau=(\tau_i)).
 \ee
Contracting both sides of (\ref{y51}) by $u^j$ we get
 $$
\tau_i=a_{ir}u^r-a_{rm}u^ru^mu_i.
 $$
Plugging the above $\tau_i$ into (\ref{y51}) we have
 $$
a_{ij}=a_{jr}u^ru_i+a_{ir}u^ru_j-a_{rm}u^ru^mu_iu_j=\big(a_{jr}u^r-\frac{1}{2}a_{rm}u^ru^mu_j\big)u_i+(i/j),
 $$
which gives (\ref{y50}). \qed

\

\noindent {\it Proof of Theorem \ref{th61} :}  Let $(M,F)$ be an
$n$-dimensional Finsler manifold. By the definition of a geodesic
circle and Proposition \ref{prop31}, we know that for any two
vectors $u,v\in T_xM$ with
 $F(u)=1$ and $g_u(u,v)=0$, there is  a unique geodesic circle
 $\gamma=\gamma(s)$ satisfying $F(\dot{\gamma}(s))=1, \gamma(0)=u$
 and $D^*_{\dot{\gamma}(0)}\dot{\gamma}=v$. Now a vector field $V$
 is concircular iff. $V$ satisfies (\ref{y29}) for any curve
 $\gamma=\gamma(s)$ with $s$ being the arc-length.

We only need to prove  (\ref{y35}) under the assumption that $V$
is a concircular vector field. Then at an arbitrarily
 fixed point $x\in M$ and unit vector $u:=\dot{\gamma}(0)\in T_xM$,  by Lemma
 \ref{lem64} together with Proposition \ref{prop31}, we have
 (\ref{y45}) for arbitrary $X\in U^{\bot}$, where $U^{\bot}$ is an $(n-1)$-dimensional space perpendicular to
$u$ under the inner product $g_u$. So (\ref{y45}) is considered as
a polynomial of $X\in U^{\bot}$ and it is equivalent to
 \be\label{y050}
\widetilde{T}^k_{ij}X^iX^j=0,\ \ \ \widetilde{Z}^k_iX^i=0,\ \ \ \
\widetilde{S}^k=0.
 \ee
Note that $X$ does no belong to the total space $T_xM$ and so
generally we don't have
$\widetilde{T}^k_{ij}=0,\widetilde{Z}^k_i=0$ from (\ref{y050}).
Here  we will use  Lemma \ref{lem65}.

For the first equation in (\ref{y050}), by Lemma \ref{lem65}, we
have
 \be\label{y52}
 \widetilde{T}^k_{ij}=\widetilde{\theta}^k_i\dot{\gamma}_j+\widetilde{\theta}^k_j\dot{\gamma}_i,
 \ \ \ \ \  \widetilde{\theta}^k_i= \widetilde{T}^k_{ir}\dot{\gamma}^r
 -\frac{1}{2}\widetilde{T}^k_{rm}\dot{\gamma}^r\dot{\gamma}^m\dot{\gamma}_i,\
 \ \  (\dot{\gamma}_i:=g_{ir}(\dot{\gamma})\dot{\gamma}^r).
 \ee
For the second equation in (\ref{y050}),  we have
 \be\label{y53}
\widetilde{Z}^k_i=\widetilde{\lambda}^k\dot{\gamma}_i,\ \ \ \ \
\widetilde{\lambda}^k=\widetilde{Z}^k_r\dot{\gamma}^r.
 \ee
Finally, we can obtain (\ref{y35}) by rewriting (\ref{y52}),
$\widetilde{S}^k=0$ and (\ref{y53}) as equations on the tangent
bundle $TM$, using the expressions of
$\widetilde{T}^k_{ij},\widetilde{Z}^k_i,\widetilde{S}^k$ in Lemma
\ref{lem64}.  In the rewriting, we should note that
 $$\mathcal{L}_VC^k_{ij}-(\mathcal{L}_{V^c}\ln F)C^k_{ij}=F\mathcal{L}_{V^c}C^k_{ij},\ \ \
(\mathcal{L}_Vg_{ij})\dot{\gamma}^i\dot{\gamma}^j=F^{-2}V^c(F^2)
$$
This completes the proof of Theorem \ref{th61}.   \qed

\subsection{Proofs of Theorems \ref{th1} and \ref{th2}}

Using Theorem \ref{th61}, we can complete the proofs of Theorems
\ref{th1} and \ref{th2}.

\

\noindent {\it Proof of Theorem \ref{th1} :} If $V$ is conformal
satisfying $\mathcal{L}_{V^c}g_{ij}=2\rho g_{ij}$ with
$\rho=\rho(x)$ being a scalar function on $M$, then we have
 $$
 0=\mathcal{L}_{V^c}\delta^i_j=\mathcal{L}_{V^c}(g^{ir}g_{rj})=(\mathcal{L}_{V^c}g^{ir})g_{rj}
 +g^{ir}(\mathcal{L}_{V^c}g_{rj})=(\mathcal{L}_{V^c}g^{ir})g_{rj}
 +2\rho\delta^i_j,
 $$
which gives $\mathcal{L}_{V^c}g^{ij}=-2\rho g^{ij}$. Thus we have
(by $\dot{\pa}_r\mathcal{L}_{V^c}=\mathcal{L}_{V^c}\dot{\pa}_r$)
 $$
\mathcal{L}_{V^c}C^i_{jk}=\mathcal{L}_{V^c}(g^{ir}C_{rjk})=(\mathcal{L}_{V^c}g^{ir})C_{rjk}+
\frac{1}{2}g^{ir}\dot{\pa}_r\mathcal{L}_{V^c}g_{jk}=0.
$$
So the Lie derivative of the mean Cartan torsion $I_i$ along $V^c$
vanishes  ($\mathcal{L}_{V^c}I_i=0$).

Conversely, first assume $V$ is concircular.    Then $V$ satisfies
(\ref{y35}) in Theorem \ref{th61}. For the first equation of
(\ref{y35}), the contraction over the indices $j$ and $k$ gives
 \be\label{y55}
 \mathcal{L}_{V^c}y_i=F^{-2}V^c(F^2)y_i+\frac{2}{3n}F^2\mathcal{L}_{V^c}I_i,
 \ee
where $n$ is the dimension of the Finsler manifold $(M,F)$. Futher
assume   the Lie derivative of the mean Cartan torsion  along
$V^c$ vanishes. Then by (\ref{y55}) we have
 \be\label{y56}
\mathcal{L}_{V^c}y_i=\lambda y_i,\ \ \ \ (\lambda:=F^{-2}V^c(F^2).
 \ee
Differentiating (\ref{y56}) by $y^j$ we obtain
 \be\label{y57}
\mathcal{L}_{V^c}g_{ij}=\lambda_{\cdot j}y_i+\lambda g_{ij},
 \ee
from which we get $\lambda_{\cdot j}y^i=\lambda_{\cdot i}y^j$.
Using this and the contraction on both sides by $y^j$, we
immediately have $\lambda_{\cdot i}=0$ since the zero homogeneity
of $\lambda$ gives $\lambda_{\cdot j}y^j=0$. So by (\ref{y57}),
$V$ is conformal satisfying $\mathcal{L}_{V^c}g_{ij}=\lambda
g_{ij}$ for a scalar function $\lambda=\lambda(x)$.    \qed

\

\noindent {\it Proof of Theorem \ref{th2} :} By assumption, $V$ is
conformal satisfying $\mathcal{L}_{V^c}g_{ij}=2\rho g_{ij}$ with
$\rho=\rho(x)$ being a scalar function on $M$. Then we have
 \beq
 &&\mathcal{L}_{V^c}C^k_{ij}=0,\ \ \
 \mathcal{L}_{V^c}y_i=2\rho y_i,\ \ \  (\mathcal{L}_{V^c}y_i)_{|j}=2\rho_j
 y_i,\ \  (\rho_i:=\rho_{x^i}),
\label{y58} \\
 &&V^c(F^2)=2\rho F^2,\ \  \   [V^c(F^2)]_{|i}=2F^2\rho_i,\ \  \
 [V^c(F^2)]_{|i|j}=2F^2\rho_{i|j}.\label{y59}
 \eeq
By $\mathcal{L}_{V^c}g_{ij}=2\rho g_{ij}$, we have (by Lemma
\ref{lem41})
 \be\label{y60}
 V_{i|0}+V_{0|i}=2\rho y_i,\ \ \   V_{i|0|0}+V_{0|i|0}=2\rho_0
 y_i,\ \ \  V_{0|0}=\rho F^2.
 \ee
Then by a Ricci identity of Cartan connection and (\ref{y60}) we
have
 \be\label{y61}
 V_{0|i|0}=V_{0|0|i}+V_{0\cdot m}R^m_{\ i}=V_{0|0|i}+V_mR^m_{\ i}=F^2\rho_i+V_mR^m_{\ i}.
 \ee
Plugging (\ref{y61}) into the second formula of (\ref{y60}) we
obtain
 \be\label{y62}
V_{i|0|0}+V_mR^m_{\ i}=2\rho_0 y_i-F^2\rho_i,\ \ \ V^i_{\
|0|0}+V^mR^i_{\ m}=2\rho_0 y^i-F^2\rho^i
 \ee
From (\ref{y62}), Lemmas \ref{lem41} and \ref{lem42}, we obtain
(note that $A^k_{j0}=B^k_{j0}$)
 \beq
 &&A^k_{00}=2\rho_0 y^k-F^2\rho^k,\ \
 A^k_{00|0}=2\rho_{0|0}y^k-F^2\rho^k_{\ |0},\label{y63}\\
 &&A^k_{j0}=\frac{1}{2}A^k_{00\cdot
 j}=\rho_jy^k+\rho_0\delta^k_j-y_j\rho^k+F^2\rho^rC^k_{jr}.\label{y64}
 \eeq

Now by Theorem \ref{th61}, we see that $V$ is concircular iff.
(\ref{y35}) holds. So we only need to simplify (\ref{y35}) with
the help of (\ref{y58}), (\ref{y59}), (\ref{y63}) and (\ref{y64}).
By (\ref{y58}) and  (\ref{y59}), we see the first equation of
(\ref{y35}) automatically holds. From (\ref{y59}) and (\ref{y63}),
the second equation of (\ref{y35}) is reduced to
 \be\label{y65}
 \rho_{0|0}y^k=F^2\rho^k_{\ |0}.
 \ee
By (\ref{y58}), (\ref{y59}), (\ref{y63}) and (\ref{y64}), the
third equation of (\ref{y35}) is reduced to $\rho^rC^k_{ir}=0$,
since we have
 $$
 Z^k_i=2F^2(F^2\rho^rC^k_{ir}-3\rho^ky_i),\ \ \
 \lambda^k=-6F^4\rho^k.
 $$
By (\ref{y65}) we have
 \be\label{y66}
 \rho^i_{\ |0}=\tau y^i, \ \ \  \rho_{i |0}=\tau y_i, \ \ \
 (\tau:=F^{-2}\rho_{0|0}).
 \ee
Differentiating (\ref{y66}) by $y^j$ we obtain (by a Ricci
identity of Berwald connection)
 \be\label{y67}
 \rho_{i;j}=\tau g_{ij}+\tau_{\cdot j}y_i\ \    (\Longleftrightarrow\
 \rho_{i|j}-\rho_rC_{ij|0}^r=\tau g_{ij}+\tau_{\cdot j}y_i),
 \ee
from which we again get $\tau_{\cdot j}y_i=\tau_{\cdot i}y_j$.
Thus we have $\tau_{\cdot i}=0$, which means that $\tau$ is a
scalar function on $M$. From (\ref{y67}) we have $\rho_{i;j}=\tau
g_{ij}$, or equivalently $\rho_{i|j}=\tau g_{ij}$ since
$\rho^rC^k_{ir}=0$ and (\ref{y66}) imply $\rho_rC_{ij|0}^r=0$. Now
we have obtained (\ref{y1}). \qed

\subsection{Proof of Theorem \ref{th3}}

We first show the following lemma which is needed in the proof of
Theorem \ref{th3} (i).

\begin{lem}\label{lem66}
Let $\widetilde{F}$ and $F$ be two conformally related Finsler
metrics on a same
 manifold $M$ with $\widetilde{F}=u^{-1}F$, and
$\gamma=\gamma(s)$ be a curve with $F(\dot{\gamma}(s))=1$. Then we
have
 \be\label{y68}
\widetilde{D}^*_{\dot{\gamma}}\widetilde{D}^*_{\dot{\gamma}}\dot{\gamma}
-3\frac{\widetilde{g}_{\dot{\gamma}}(\dot{\gamma},\widetilde{D}^*_{\dot{\gamma}}\dot{\gamma})}
{\widetilde{F}^2(\dot{\gamma})}\
\widetilde{D}^*_{\dot{\gamma}}\dot{\gamma}=D^*_{\dot{\gamma}}D^*_{\dot{\gamma}}\dot{\gamma}
+\frac{1}{u}D^*_{\dot{\gamma}}U+\lambda \dot{\gamma},
 \ee
 where $U=u^i\pa_i$ is a vector field along $\gamma$ defined by $u^i:=g^{ir}(\dot{\gamma})u_r$
 with $u_r:=u_{x^r}$, and $\lambda=\lambda(s)$ is a function along
 $\gamma$.
\end{lem}

{\it Proof :} Since $\widetilde{F}=u^{-1}F$, a direct computation
from (\ref{y5}) shows that
 \be\label{y69}
 \widetilde{G}^i=G^i-\frac{1}{u}u_0y^i+\frac{1}{2u}F^2u^i,\ \ \
  \widetilde{G}^i_j=G^i_j-\frac{1}{u}\big(u_jy^i+u_0\delta^i_j-y_ju^i+F^2C^i_{jr}u^r\big).
 \ee
By the first formula of (\ref{y69}) we have
 \be\label{y70}
(\widetilde{D}^*_{\dot{\gamma}}\dot{\gamma})^k=\ddot{\gamma}^k+2\widetilde{G}^k
=(D^*_{\dot{\gamma}}\dot{\gamma})^k-\frac{2}{u}g_{\dot{\gamma}}(\dot{\gamma},U)\dot{\gamma}^k+\frac{1}{u}u^k.
 \ee
Then by (\ref{y70}), we first have
 \be\label{y71}
\widetilde{D}^*_{\dot{\gamma}}\widetilde{D}^*_{\dot{\gamma}}\dot{\gamma}
=\widetilde{D}^*_{\dot{\gamma}}D^*_{\dot{\gamma}}\dot{\gamma}
-\frac{d}{ds}\big[\frac{2}{u}g_{\dot{\gamma}}(\dot{\gamma},U)\big]\dot{\gamma}
-\frac{2}{u}g_{\dot{\gamma}}(\dot{\gamma},U)\widetilde{D}^*_{\dot{\gamma}}\dot{\gamma}+\frac{d}{ds}(\frac{1}{u})U
+\frac{1}{u}\widetilde{D}^*_{\dot{\gamma}}U.
 \ee
We respectively have
 \beq
(\widetilde{D}^*_{\dot{\gamma}}D^*_{\dot{\gamma}}\dot{\gamma})^k
\hspace{-0.6cm}&&=\frac{d}{ds}(D^*_{\dot{\gamma}}\dot{\gamma})^k
+\dot{\gamma}^j(D^*_{\dot{\gamma}}\dot{\gamma})^i\widetilde{\Gamma}_{ij}^{*k}=\frac{d}{ds}(D^*_{\dot{\gamma}}\dot{\gamma})^k
+\dot{\gamma}^j(D^*_{\dot{\gamma}}\dot{\gamma})^i(^*\widetilde{\Gamma}_{ij}^k+\widetilde{C}^k_{ir}\widetilde{Y}^r_j)\nonumber\\
&&=\frac{d}{ds}(D^*_{\dot{\gamma}}\dot{\gamma})^k
+(D^*_{\dot{\gamma}}\dot{\gamma})^i\big[\widetilde{G}^k_i+\widetilde{C}^k_{ir}(\widetilde{D}^*_{\dot{\gamma}}\dot{\gamma})^r\big],\label{y72}\\
(\widetilde{D}^*_{\dot{\gamma}}U)^k\hspace{-0.6cm}&&=\frac{d}{ds}u^k
+u^i\big[\widetilde{G}^k_i+\widetilde{C}^k_{ir}(\widetilde{D}^*_{\dot{\gamma}}\dot{\gamma})^r\big],\label{y73}
 \eeq
By $\widetilde{C}^i_{jk}=C^i_{jk}$ (conformally invariant), the
second formula of (\ref{y69}) and (\ref{y70}), we can rewrite
(\ref{y72}) and (\ref{y73}) as follows:
 \beq
\widetilde{D}^*_{\dot{\gamma}}D^*_{\dot{\gamma}}\dot{\gamma}
\hspace{-0.6cm}&&=D^*_{\dot{\gamma}}D^*_{\dot{\gamma}}\dot{\gamma}
-\frac{1}{u}\big\{g_{\dot{\gamma}}(U,D^*_{\dot{\gamma}}\dot{\gamma})\dot{\gamma}
+g_{\dot{\gamma}}(U,\dot{\gamma})D^*_{\dot{\gamma}}\dot{\gamma}\big\},\label{y74}\\
\widetilde{D}^*_{\dot{\gamma}}U\hspace{-0.6cm}&&=D^*_{\dot{\gamma}}U-\frac{1}{u}g_{\dot{\gamma}}(U,U)\dot{\gamma}.\label{y75}
 \eeq
Plugging (\ref{y70}),  (\ref{y74}) and (\ref{y75}) into
(\ref{y71}), we obtain
 \be\label{y76}
\widetilde{D}^*_{\dot{\gamma}}\widetilde{D}^*_{\dot{\gamma}}\dot{\gamma}
=D^*_{\dot{\gamma}}D^*_{\dot{\gamma}}\dot{\gamma}+\frac{1}{u}D^*_{\dot{\gamma}}U
-\frac{3}{u}g_{\dot{\gamma}}(\dot{\gamma},U)\big(\widetilde{D}^*_{\dot{\gamma}}\dot{\gamma}+\frac{1}{u}U\big)
+\lambda_1\dot{\gamma},
 \ee
where $\lambda_1=\lambda_1(s)$ is a function along $\gamma$. By
$\widetilde{F}=u^{-1}F$ and (\ref{y70}), it is easy to see that
 \be\label{y77}
-3\frac{\widetilde{g}_{\dot{\gamma}}(\dot{\gamma},\widetilde{D}^*_{\dot{\gamma}}\dot{\gamma})}
{\widetilde{F}^2(\dot{\gamma})}\
\widetilde{D}^*_{\dot{\gamma}}\dot{\gamma}=
\frac{3}{u}g_{\dot{\gamma}}(\dot{\gamma},U)\big(\widetilde{D}^*_{\dot{\gamma}}\dot{\gamma}+\frac{1}{u}U\big)
+\lambda_2\dot{\gamma},
 \ee
where $\lambda_2=\lambda_2(s)$ is a function along $\gamma$. By
(\ref{y76}) and (\ref{y77}), we immediately obtain (\ref{y68}).
This completes the proof.   \qed

\

Now we can get started with the proof of Theorem \ref{th3} (i).
Let $\widetilde{F}$ be  conformally related to $F$ satisfying
$\widetilde{F}=u^{-1}F$ on a same manifold $M$.

Assume $u$ satisfies (\ref{y3}). Let $\gamma=\gamma(s)$ be an
arbitrary geodesic circle of $(M,F)$ with $F(\dot{\gamma}(s))=1$.
Then by the definition of a geodesic circle, we see
$D^*_{\dot{\gamma}}D^*_{\dot{\gamma}}\dot{\gamma}$ is parallel to
$\dot{\gamma}$. By (\ref{y7}) and then by (\ref{y3}), we have
 $$
 (D^*_{\dot{\gamma}}U)^k=\dot{\gamma}^ju^k_{/j}=\dot{\gamma}^j(u^k_{|j}+u^k|_rY^r_j)
 =\dot{\gamma}^ju^k_{|j}+u^mC^k_{mr}(D^*_{\dot{\gamma}}\dot{\gamma})^r=\lambda
 \dot{\gamma}^k.
 $$
 By Lemma
\ref{lem66} we have (\ref{y68}). Now it is easy to see that
(\ref{y68}) implies that the following vector
 $$
\widetilde{D}^*_{\dot{\gamma}}\widetilde{D}^*_{\dot{\gamma}}\dot{\gamma}
-3\frac{\widetilde{g}_{\dot{\gamma}}(\dot{\gamma},\widetilde{D}^*_{\dot{\gamma}}\dot{\gamma})}
{\widetilde{F}^2(\dot{\gamma})}\
\widetilde{D}^*_{\dot{\gamma}}\dot{\gamma}
 $$
is parallel to $\dot{\gamma}$. Thus by Proposition \ref{prop37},
the curve $\gamma$ is also a geodesic circle (as points set) of
$(M,\widetilde{F})$. Similarly, a geodesic circle of
$\widetilde{F}$ is also a geodesic circle of $F$. This means that
$\widetilde{F}$ and $F$ are concircular.

Conversely, suppose that $\widetilde{F}$ and $F$ are concircular.
Then an arbitrary geodesic circle $\gamma$ of $(M,F)$ is also a
geodesic circle of $(M,\widetilde{F})$. So by Proposition
\ref{prop37} and (\ref{y68}), we see that $D^*_{\dot{\gamma}}U$ is
parallel to $\dot{\gamma}$, which shows that
 \be\label{y78}
\big((D^*_{\dot{\gamma}}U)^k=\dot{\gamma}^ju^k_{/j}
 =\big) \  \dot{\gamma}^ju^k_{|j}+u^mC^k_{mr}(D^*_{\dot{\gamma}}\dot{\gamma})^r=\lambda
 \dot{\gamma}^k,
 \ee
where $\lambda=\lambda(\gamma,\dot{\gamma})$ is a scalar function
along $\gamma$, and actually the contraction of (\ref{y78}) by
$\dot{\gamma}_k:=g_{kr}(\dot{\gamma})\dot{\gamma}^r$ gives
$\lambda=\dot{\gamma}^i\dot{\gamma}^ju_{i|j}$. Let
$w:=\dot{\gamma}$ and $W^{\bot}$ be the $(n-1)$-dimensional space
perpendicular to $w$ with respect to the inner product $g_w$. Then
for fixed $\dot{\gamma}$, by Proposition \ref{prop31}, we see that
(\ref{y78}) is a polynomial equation of the variable
$X:=D^*_{\dot{\gamma}}\dot{\gamma}\in W^{\bot}$. Thus (\ref{y78})
is equivalent to
 \be\label{y79}
  \dot{\gamma}^ju^k_{|j}=\lambda
\dot{\gamma}^k,\ \ \ \ \ \ \
u^mC^k_{mr}(D^*_{\dot{\gamma}}\dot{\gamma})^r=0.
 \ee
 We can write (\ref{y79}) as equations on the tangent bundle $TM$
as follows:
 \be\label{y80}
 u^k_{\ |0}=\lambda y^k,\ \ \ \  \ \ \   u^mC^k_{mr}=\tau^ky_r.
 \ee
The first equation in (\ref{y80}) is similar to (\ref{y66}). So
$\lambda=\lambda(x)$ is a scalar function on $M$, and then
$u_{i;j}=\lambda g_{ij}$. For the second equation of (\ref{y80}),
the contraction by $y^r$ immediately gives $\tau^k=0$ and thus
$u^mC^k_{mr}=0$. Now we have proved (\ref{y3}).

\

Before the proof of Theorem \ref{th3} (ii), we first give a brief
introduction for some basic points needed here. It is known that
 if two sprays $\widetilde{G}^i$ and
$G^i$ satisfy $\widetilde{G}^i=G^i+H^i$, then their Riemann
curvature tensors $\widetilde{R}^i_{\ k}$ and $R^i_{\ k}$ are
related by
 \be\label{y81}
 \widetilde{R}^i_{\ k}=R^i_{\ k}+2H^i_{\ ;k}-y^mH^i_{\ ;m\cdot k}+2H^mH^i_{\ \cdot m\cdot k}-H^i_{\ \cdot m}H^m_{\ \cdot k},
 \ee
where the symbol $_;$ denotes the horizontal covariant derivative
of Berwald connection of $G^i$. A Finsler metric $F$ is said to be
of scalar (resp. isotropic) flag curvature, if the Riemann
curvature satisfies
 \be\label{y87}
R^i_{\ k}=K(F^2\delta^i_k-y^iy_k),
 \ee
where $K=K(x,y)$ is a scalar function on $TM$ (resp. $K=K(x)$ is a
scalar function on $M$). If $K$ in (\ref{y87}) is a constant, then
$F$ is called of constant flag curvature. A Finsler metric $F$ is
said to be  an Einstein metric, if the Ricci curvature is of
isotropic Ricci scalar in the following form,
 \be\label{y087}
 Ric=(n-1)KF^2,
 \ee
where $K=K(x)$  is a scalar function on $M$.

Now we show the proof. Since $\widetilde{F}=u^{-1}F$, the sprays
$\widetilde{G}^i$ and $G^i$ are related by (\ref{y69}) with $H^i$
being given by
 \be\label{y82}
H^i=-\frac{1}{u}u_0y^i+\frac{1}{2u}F^2u^i.
 \ee
Plugging (\ref{y82}) into (\ref{y81}), we obtain by a direct
computation
 \beq
\widetilde{R}^i_{\ k}\hspace{-0.6cm}&&=R^i_{\
k}+\frac{uu_{0;0}-(u_mu^m)F^2}{u^2}\delta^i_k+\frac{1}{u}F^2u^i_{
;k}+\frac{u_mu^m}{u^2}y^iy_k-\frac{1}{u}(y^iu_{k;0}+y_ku^i_{;0})\nonumber\\
&&-\frac{u^mu^r}{u^2}F^2(y^iC_{kmr}+y_kC^i_{mr})+\frac{1}{u^2}F^2(uu^r_{;0}-3u_0u^r)C^i_{kr}+\frac{1}{u}F^2u^rC^i_{kr;0}\nonumber\\
&&+\frac{u^ru^m}{u^2}F^4(C^i_{pr}C^p_{km}-C^i_{mr\cdot
k}).\label{y83}
 \eeq
Then by (\ref{y83}), the Ricci curvatures
$\widetilde{Ric}:=\widetilde{R}^m_{\ m}$ and $Ric:=R^m_{\ m}$ are
related by
 \beq\label{y84}
\widetilde{Ric}\hspace{-0.6cm}&&=Ric+\frac{n-2}{u}u_{0;0}+\frac{1}{u^2}\big[uu^m_{;m}-(n-1)u^mu_m+uI^ru_{r;0}+u^r(uI_{r;0}-3u_0I_r)\big]F^2\nonumber\\
&&-\frac{1}{u^2}u^ru^m(C^i_{jm}C^j_{ir}-2I^iC_{imr}+I_{m\cdot
r})F^4
 \eeq

Now suppose $\widetilde{F}$ and $F$ are concircular. Then by
Theorem \ref{th3} (i), we have (\ref{y3}). Plugging (\ref{y3})
into (\ref{y83}) and (\ref{y84}), we respectively have
 \beq
 \widetilde{R}^i_{\ k}\hspace{-0.6cm}&&=R^i_{\ k}+u^{-2}(2\lambda
 u-u_mu^m)(F^2\delta^i_k-y^iy_k),\label{y85}\\
  \widetilde{Ric}\hspace{-0.6cm}&&=Ric+(n-1)u^{-2}(2\lambda
 u-u_mu^m)F^2.\label{y86}
 \eeq
If $F$ is of scalar (resp. isotropic) flag curvature satisfying
(\ref{y87}), or an Einstein metric satisfying (\ref{y087}), then
plugging (\ref{y87}) into (\ref{y85}), and (\ref{y087}) into
(\ref{y86}), respectively we obtain
 \beq
\widetilde{R}^i_{\ k}\hspace{-0.6cm}&&=(Ku^2+2\lambda
 u-u_mu^m)(\widetilde{F}^2\delta^i_k-y^i\widetilde{y}_k),\label{y89}\\
 \widetilde{Ric}\hspace{-0.6cm}&&=(n-1)(Ku^2+2\lambda
 u-u_mu^m)\widetilde{F}^2.\label{y90}
 \eeq
Note that we have $u^i_{\cdot k}=0$ from the second equation in
(\ref{y3}). Now it follows from (\ref{y89}) that $\widetilde{F}$
is of scalar (resp. isotropic) flag curvature $\widetilde{K}$
given by (\ref{y4}), or from (\ref{y90}) that $\widetilde{F}$ is
an Einstein metric with the Ricci scalar $\widetilde{K}$ given by
(\ref{y4}).    \qed

\begin{rem}\label{rem61}
 In Theorem \ref{th3} (i), if $F$ is locally Euclidean, then we can solve (\ref{y3}) in a local coordinate
 such that $\widetilde{F}$ is locally expressed as $\widetilde{F}=u^{-1}|y|$.
 So $u_{i|j}=\lambda g_{ij}$ is equivalent to
 $u_{x^ix^j}=\lambda\delta_{ij}$. By integrability, we see
 $\lambda$ is a constant, and thus we obtain
  \be\label{y91}
\widetilde{ F}=\big(a|x|^2+\langle b,x\rangle+c\big)^{-1}|y|,\ \ \
 (a:=\lambda/2),
  \ee
  where $a,c$ are constant numbers and $b$ is a constant
  $n$-vector such that $u>0$.

  For convenience, suppose (\ref{y91}) is defined on the whole
  $R^n$. Let $\gamma=\xi s+\tau$ be a geodesic in the Euclidean space
  $(R^n,F)$, where $\xi,\tau$ are $n$-vectors satisfying
  $|\xi|=1$. Let $t$ be the arc-length of $\gamma$ with respect to
  $\widetilde{F}$. Then a direct computation from (\ref{y69})
  gives
   $$
\widetilde{D}^*_{\gamma'(t)}\gamma'(t)=u\big[2a\tau+b-\langle
2a\tau+b,\xi\rangle \xi\big].
   $$
So $\gamma$ is also a geodesic of $\widetilde{F}$ iff. $2a\tau+b$
is tangent to $\gamma$. Otherwise, $\gamma$ is a circle of
$\widetilde{F}$.

Let $\gamma=\gamma(s)$ be a circle of $F$. Then by Example
\ref{ex36}, $\gamma$ is written as
$$\gamma=\xi \cos ks+\eta \sin
ks+\tau,\ \ \ \ |\xi|=|\eta|=1/k).$$ Similarly, a direct
computation gives
 \be\label{y92}
\widetilde{D}^*_{\gamma'(t)}\gamma'(t)=u\big[(A\cos ks-k^2\langle
B,\xi\rangle)\xi+(A\sin ks-k^2\langle B,\eta\rangle)\eta+B\big],
 \ee
 where $t$ is the arc-length of $\gamma$ with respect to
  $\widetilde{F}$, and $A,B$ are defined by
  $$
 A:=a-k^2(\langle b,\tau\rangle+a|\tau|^2+c),\ \ \ \ B:=b+2a\tau.
  $$
  By (\ref{y92}), we easily obtain
   $$
 |\widetilde{D}^*_{\gamma'(t)}\gamma'(t)|^2_{\widetilde{g}_{\gamma'(t)}}=-\big(\langle
 B,\xi\rangle^2+\langle
 B,\eta\rangle^2\big)k^2+|B|^2+\frac{A^2}{k^2}.
   $$
   Thus we can determine the conditions for $\gamma$ to be
   a geodesic or a circle of $\widetilde{F}$.
 \end{rem}

\section{Some examples}\label{sec7}

 In this section, we give some examples to show that concircular
 vector fields might not be conformal and conformal vector fields
 might not be concircular.

\begin{ex}\label{ex71}
 Let $F=\alpha+\beta$ be an $n$-dimensional Randers metric. By the first equation in (\ref{y35}), we can prove that if $n\ge 3$, then
 any concircular vector field of $F$ must be conformal. While in
 dimension $n=2$, there exist non-conformal concircular vector
 fields, which will be exemplified as follows.

  Define   a two-dimensional Minkowskin Randers metric $F=\alpha+\beta$ by
 $$
 \alpha:=\sqrt{(y^1)^2+(y^2)^2},\ \ \ \beta:=by^1, \ \ \  (a \ constant \ b \ with \
 0<b<1),
 $$
 and a vector field $V=(V^1,V^2)$  by
  $$
 V^1:=qx^2+\eta^1,\ \ \  V^2:=-qx^1+\eta^2,
  $$
  where $q$ is a non-zero constant and $\eta=(\eta^1,\eta^2)$ is a
  constant vector. It can be easily checked that
   $$
 V^c(\alpha^2)=0,\  \ \ V^c(\beta)=bqy^2\ne0,
   $$
   which implies that $V$ is not conformal in $F$ (cf. \cite{Y3,Y4}). On the other hand, a direct verification shows that
   $V$ satisfies (\ref{y35}), and thus $V$ is
   concircular in $F$ by Theorem \ref{th61}.
\end{ex}

\begin{ex}\label{ex72}
  Let $F=e^{\sigma(x)/2}|y|$ be an $n(\ge3)$-dimensional
conformally flat Riemann metric and $V$ be a conformal vector
field of $F$ with the conformal factor $\rho=\rho(x)$. Then
 $$
 V^i=-2\big(\lambda+\langle
 d,x\rangle\big)x^i+|x|^2d^i+q_r^ix^r+\eta^i,\ \
\rho=-2(\lambda+\langle
 d,x\rangle)+\frac{1}{2}V(\sigma),
 $$
where $\lambda$ is a constant number, $d,\eta$ are constant
vectors and
  $(q_i^j)$ is skew-symmetric (cf. \cite{Y1,Y2,Y3,Y4}). If $F$ is of constant sectional curvature $\mu$ ($
 \sigma=ln4/(1+\mu|x|^2)^2$), then $V$ is also concircular by (\ref{y1}) (cf. \cite{Ya1}). Taking
 $\sigma=|x|^2$ (or many other functions), we can check that $V$ is non-concircular by (\ref{y1}).
\end{ex}

\begin{ex}\label{ex73}
 Define a projectively flat Randers metric $F=\alpha+\beta$ and a vector field $V$
 by
  \beqn
 \alpha:&&\hspace{-0.6cm}=\frac{2}{1+\mu|x|^2}|y|, \ \ \ \ \ \ \ \
 \beta:=\frac{1}{\lambda(1-\mu|x|^2)+\langle
 d,x\rangle}\Big\{\langle d,y\rangle-\frac{2\mu(2\lambda+\langle d,x\rangle)\langle x,
 y\rangle}{1+\mu|x|^2}\Big\},\\
 V^i:&&\hspace{-0.6cm}=-2\big(\lambda+\langle
 d,x\rangle\big)x^i+|x|^2d^i,
  \eeqn
  where the constant $\lambda$ and the constant vector $d=(d^i)\ne 0$
  satisfy $|d|^2+4\mu\lambda^2=0$. It has been verified in
  \cite{Y3,Y4} that
 $V$ is a non-homothetic conformal vector field in $F$ with the conformal factor $\rho$ given by
   $$\rho=-\frac{2\big[\lambda(1-\mu|x|^2)+\langle d,x\rangle\big]}{1+\mu|x|^2}.$$
   It can be checked directly that $\rho$ does not satisfy
   (\ref{y1}), and so $V$ is not concircular in $F$ by Theorem
   \ref{th2}.
\end{ex}

\noindent Zhongmin Shen\\
Department of Mathematical Sciences\\
Indinan University Purdue University Indianapolis (IUPUI)\\
402 N. Blackford Street\\
Indianapolis, IN 46202-3216, USA\\
zshen@math.iupui.edu

\vspace{0.5cm}

\noindent Guojun Yang \\
Department of Mathematics \\
Sichuan University \\
Chengdu 610064, P. R. China \\
yangguojun@scu.edu.cn


\begin{thebibliography}{999}


\bibitem{AIM} P. L. Antonelli, R. S. Ingarden and M. Matsumoto, {\it The Theory of
Sprays and Finsler Spaces with Applications in Physics and
Biology}, KLuwer Academic Publishers, London, 1993.

\bibitem{BS} B. Bidabad and Z. Shen, Circle-preserving transformations in
Finsler spaces, Publ. Math. Debrecen, {\bf 81} (2012), 435-445.

\bibitem{De} A. Deicke, Uber die Finsler-Raume mit $A_i = 0$, Arch.
Math., {\bf 4} (1953), 45-51.

\bibitem{Fe} J. Ferrand, {\it Concircular transformations of
Riemannian manifolds}, Ann. Acad. Sci. Fennice-M., {\bf 10}
(1985), 163-171.

\bibitem{Is} S. Ishihara, On infinitesimal concircular
transformations, K\^{o}dai Math. Sem. Rep. {\bf 12} (1960), 45-56.

\bibitem{JB} P. Joharinad and B. Bidabvad, Conformal vector fields
on Finsler spaces, Diff. Geom. Appl., {\bf 31} (2013), 33-40.

\bibitem{IT} S. Ishihara and Y. Tashiro, On Riemannian manifolds
admitting a concircular transformation, Math. J. Okayama Univ.
{\bf 9} (1959), 19-47.

\bibitem{Ma} M. Matsumoto, {\it Theory of Y-extremal and minimal hypersurfaces in
a Finsler space}, J. Math. Kyoto Univ., {\bf 26}(4)(1986),
647-665.

\bibitem{NY} K. Nomizu and K. Yano, On Circles and Spheres in Riemannian
Geometry, Math. Ann., {\bf 210} (1974), 163-170.

\bibitem{Ta} Y. Tashiro, Complete Riemannian manifolds and some
vector fields, Trans. Amer. math. Soc. {\bf 117} (1965), 251-275.

\bibitem{Yan} K. Yano, On circular geometry, I. Concircular
transformations. Proc. Imp. Acad. Tokyo, {\bf 16} (1940), 195-200;
II. Integrability conditions of $\rho_{\mu\lambda}=\phi
g_{\mu\lambda}$, ibid. {\bf 16} (1940), 354-360; III. Theory of
curves, ibid. {\bf 16} (1940), 442-448; IV. Theory of subspaces,
ibid. {\bf 16} (1940), 505-511; V. Einstein spaces, ibid. {\bf 18}
(1942), 446-451.



\bibitem{MH} X. Mo and L. Huang, {\it On curvature decreasing property of a
class of navigation problems}, Publ. Math. Debrecen, {\bf 71}
(1-2) (2007), 141-163.





\bibitem{TI} S. Tachibana and S. Ishihara, {\it On infininitesimal
h'olomorphically projective transformations in Kaehlerian
manifolds}, Tohoku Math. Journ. {\bf 12} (1960), 77-101.



\bibitem{Vo} W. O. K. Vogel, {\it Transformationen in Riemannschen R\"{a}umen}, (German)
Arch. Math. (Basel) 21 (1970-71), 641--645.


\bibitem{Ya1} K. Yano, {\it Einstein spaces admitting a
one-parameter group of conformal transformations}, Ann. Math.,
{\bf 69}(2) (1959), 451-460.

\bibitem{Ya} K. Yano, {\it The Theory of Lie Derivatives and its Applications}, North-Holland Pub., 1957.

\bibitem{Y1} G. Yang, {\it On Randers metrics of isotropic
S-curvature}, Acta Math. Sin., {\bf 52}(6)(2009), 1147-1156 (in
Chinese).


\bibitem{Y2} G. Yang, {\it On Randers metrics of isotropic
S-curvature II}, Publ. Math. Debreceen, {\bf 78}(1) (2011), 71-87.


\bibitem{Y3} G. Yang, {\it Conformal Vector Fields On Projectively Flat $(\alpha,\beta)$-Finsler Spaces}, preprint.

\bibitem{Y4} G. Yang, {\it Conformal Vector Fields of a Class of  Finsler
Spaces}, preprint.






\end{thebibliography}
\end{document}